\documentclass{amsart}
\usepackage[english]{babel}
\usepackage[all]{xy}
\usepackage{xkeyval}
\usepackage{graphicx}
\usepackage{amssymb,amscd,verbatim, amsthm,amsmath,amsgen,xspace}
\usepackage{latexsym}
\usepackage{amsfonts}
\usepackage{color}
\usepackage{ulem}
\usepackage[all]{xy}
\usepackage{latexsym}
\usepackage{exscale}
\usepackage{comment}

\usepackage{ytableau}
\usepackage{color}

\usepackage[all]{xy}

\usepackage{amsbsy}

\def\Bbb{\mathbb}

\newtheorem{prop}[equation]{Proposition}
\newtheorem*{prop*}{Proposition}
\newtheorem{thm}[equation]{Theorem}
\newtheorem*{thm*}{Theorem}
\newtheorem{lem}[equation]{Lemma}
\newtheorem*{lem*}{Lemma}

\newtheorem*{kor*}{Corollary}
\newtheorem{cor}[equation]{Corollary}
\newtheorem{rem}[equation]{Remark}
\newtheorem{df}[equation]{Definition}
\newtheorem{ex}[equation]{Example}

\numberwithin{equation}{section}

\newcommand{\Id}{\operatorname{id}}

\newcommand{\frg}{\mathfrak{g}}
\newcommand{\frh}{\mathfrak{h}}
\newcommand{\frk}{\mathfrak{k}}

\newcommand{\frn}{\mathfrak{n}}

\newcommand{\frp}{\mathfrak{p}}

\newcommand{\frt}{\mathfrak{t}}
\newcommand{\fru}{\mathfrak{u}}

\newcommand{\frgl}{\mathfrak{gl}}
\newcommand{\frsl}{\mathfrak{sl}}
\newcommand{\frsu}{\mathfrak{su}}
\newcommand{\frso}{\mathfrak{so}}
\newcommand{\frsp}{\mathfrak{sp}}

\newcommand{\Cas}{\operatorname{Cas}}

\newcommand{\bbar}{\,|\,}
\newcommand{\del}{\partial}

\newcommand \bul{\bullet}

\def\bbC{\mathbb{C}}

\def\bbR{\mathbb{R}}
\def\bbZ{\mathbb{Z}}

\let\ccdot\cdot
\def\cdot{\hbox to 2.5pt{\hss$\ccdot$\hss}}

\newcommand{\bC}{{\Bbb C}}

\newcommand{\eq}{\begin{equation}}
\newcommand{\eeq}{\end{equation}}
\newcommand{\eqn}{\begin{equation*}}
\newcommand{\bmul}{\begin{multline*}}
\newcommand{\eemul}{\end{multline*}}
\newcommand{\eeqn}{\end{equation*}}
\newcommand{\pf}{\begin{proof}}
\newcommand{\epf}{\end{proof}}

\newcommand{\pa}{\partial}

\newcommand{\la}{\lambda}

\renewcommand{\phi}{\varphi}

\newcommand{\La}{\Lambda}

\newcommand{\eps}{\varepsilon}
\newcommand{\half}{\frac{1}{2}}
\newcommand{\dom}{\operatorname{dom}}

\let\ssize\scriptstyle
\newif\ifFIRST\newdimen\MAXright\MAXright0pt
\def\sdynkin{\bgroup\eightpoint\dynkin}
\def\endsdynkin{\enddynkin\egroup}
\def\dynkin{\bgroup\FIRSTtrue\hskip.5em\setbox1\hbox{$\diagup$}%
	\setbox2\hbox{$\diagdown$}%
	\setbox0\hbox to2\wd1{\hrulefill}%
	\setbox3\hbox{$\bullet$}%
	\setbox4\hbox{$\times$}%
	\setbox7\hbox{$\circ$}
	\def\whiteroot##1{\ifFIRST\setbox5\hbox{$##1$}\ifdim\wd5>1.3em
		\hskip-.5em\hskip.5\wd5\fi\fi\FIRSTfalse
		\hskip-.25em\raise1.5\wd3\hbox to0pt{\hss\hskip.45em$
			\ssize##1$\hss}\copy7\hskip-.25em\setbox6\hbox{$##1$}
		\MAXright\wd6}
	\def\root##1{\ifFIRST\setbox5\hbox{$##1$}\ifdim\wd5>1.3em%
		\hskip-.5em\hskip.5\wd5\fi\fi\FIRSTfalse%
		\hskip-.25em\raise1.5\wd3\hbox to0pt{\hss\hskip.45em$%
			\ssize##1$\hss}\copy3\hskip-.25em\setbox6\hbox{$##1$}%
		\MAXright\wd6}%
	\def\whitedroot##1{\ifFIRST\setbox5\hbox{$##1$}\ifdim\wd5>1.3em
		\hskip-.5em\hskip.5\wd5\fi\fi\FIRSTfalse
		\hskip-.25em\lower1.8\wd3\hbox to0pt{\hss\hskip.45em$
			\ssize##1$\hss}\copy7\hskip-.25em\setbox6\hbox{$##1$}
		\MAXright\wd6}%
	\def\whiterroot##1{\hskip-.25em\copy7\hbox to0pt{\hskip.3em$\ssize##1$\hss}%
		\hskip-.25em\setbox6\hbox{\hskip.6em$##1##1$}%
		\MAXright\wd6}%
	\def\droot##1{\ifFIRST\setbox5\hbox{$##1$}\ifdim\wd5>1.3em%
		\hskip-.5em\hskip.5\wd5\fi\fi\FIRSTfalse%
		\hskip-.25em\lower1.8\wd3\hbox to0pt{\hss\hskip.45em$%
			\ssize##1$\hss}\copy3\hskip-.25em\setbox6\hbox{$##1$}%
		\MAXright\wd6}%
	\def\rroot##1{\hskip-.25em\copy3\hbox to0pt{\hskip.3em$\ssize##1$\hss}%
		\hskip-.25em\setbox6\hbox{\hskip.6em$##1##1$}%
		\MAXright\wd6}%
	\def\norroot##1{\hskip-.36em\copy4\hbox to0pt{\hskip.3em$\ssize##1$\hss}%
		\hskip-.48em\setbox6\hbox{\hskip.6em$##1##1$}%
		\MAXright\wd6}%
	\def\noroot##1{\ifFIRST\setbox5\hbox{$##1$}\ifdim\wd5>1.3em%
		\hskip-.5em\hskip.5\wd5\fi\fi\FIRSTfalse%
		\hskip-.36em\raise1.5\wd3\hbox to0pt{\hss\hskip.6em$%
			\ssize##1$\hss}\copy4\hskip-.38em\setbox6\hbox{$##1$}%
		\MAXright\wd6}%
	\def\nodroot##1{\ifFIRST\setbox5\hbox{$##1$}\ifdim\wd5>1.3em%
		\hskip-.5em\hskip.5\wd5\fi\fi\FIRSTfalse%
		\hskip-.36em\lower1.8\wd3\hbox to0pt{\hss\hskip.6em$%
			\ssize##1$\hss}\copy4\hskip-.38em\setbox6\hbox{$##1$}%
		\MAXright\wd6}%
	\def\nolink{\hskip\wd0}
	\def\link{\raise.22em\copy0}%
	\def\llink##1{\raise.32em\copy0\hskip-\wd0%
		\raise.12em\copy0\hskip-.5\wd0\hbox to0pt{\hss$##1$\hss}\hskip.5\wd0}%
	\def\lllink##1{\raise.22em\copy0\hskip-\wd0\raise.32em\copy0\hskip-\wd0%
		\raise.12em\copy0\hskip-.5\wd0\hbox to0pt{\hss$##1$\hss}\hskip.5\wd0}%
	\def\rootupright##1{\hbox to0pt{\raise.45em\copy1\hskip-.25em\raise1.3\ht1%
			\hbox{\copy3\hskip.3em$\ssize##1$}\hss}%
		\setbox6\hbox{\hskip.6em\copy1\copy1$##1##1$}%
		\ifdim\MAXright<\wd6\MAXright\wd6\fi}%
	\def\whiterootupright##1{\hbox to0pt{\raise.45em\copy1\hskip-.25em\raise1.3\ht1
			\hbox{\copy7\hskip.3em$\ssize##1$}\hss}
		\setbox6\hbox{\hskip.6em\copy1\copy1$##1##1$}
		\ifdim\MAXright<\wd6\MAXright\wd6\fi}
	\def\norootupright##1{\hbox to0pt{\raise.45em\copy1\hskip-.36em\raise1.3\ht1%
			\hbox{\copy4\hskip.3em$\ssize##1$}\hss}%
		\setbox6\hbox{\hskip.6em\copy1\copy1$##1##1$}%
		\ifdim\MAXright<\wd6\MAXright\wd6\fi}%
	\def\rootdownright##1{\hbox to0pt{\raise-.5em\copy2\hskip-.25em\raise-1.35\ht1%
			\hbox{\copy3\hskip.3em$\ssize##1$}\hss}\setbox6%
		\hbox{\hskip.6em\copy2\copy2$##1##1$}%
		\ifdim\MAXright<\wd6\MAXright\wd6\fi}%
	\def\whiterootdownright##1{\hbox to0pt{\raise-.5em\copy2\hskip-.25em\raise-1.35\ht1
			\hbox{\copy7\hskip.3em$\ssize##1$}\hss}\setbox6
		\hbox{\hskip.6em\copy2\copy2$##1##1$}
		\ifdim\MAXright<\wd6\MAXright\wd6\fi}
	\def\rootdown##1{\hbox to0pt{\hskip-.05em\vrule height.25em depth.65em%
			\hskip-.25em\raise-.95em\hbox{\copy3\hskip.3em$\ssize##1$}\hss}%
		\setbox6\hbox{$##1$}%
		\ifdim\MAXright<\wd6\MAXright\wd6\fi}%
	\def\whiterootdown##1{\hbox to0pt{\hskip-.05em\vrule height.25em depth.65em
			\hskip-.25em\raise-.95em\hbox{\copy7\hskip.3em$\ssize##1$}\hss}
		\setbox6\hbox{$##1$}
		\ifdim\MAXright<\wd6\MAXright\wd6\fi}
	\def\dots{\hskip.5em\cdots\hskip.5em}}%
\def\enddynkin{\ifdim\MAXright>1em\hskip.5\MAXright\else\hskip.5em\fi\egroup}%
 
\begin{document} 

\title[Unitary highest weight modules]{On the classification of unitary highest weight modules}
\author{Pavle Pand\v zi\'c}
\address[Pand\v zi\'c]{Department of Mathematics, Faculty of Science, University of Zagreb, Bijeni\v cka 30, 10000 Zagreb, Croatia}
\email{pandzic@math.hr}
\author{Ana Prli\'c}
\address[Prli\'c]{Department of Mathematics, Faculty of Science, University of Zagreb, Bijeni\v cka 30, 10000 Zagreb, Croatia}
\email{anaprlic@math.hr}
\author{Vladim\'{\i}r Sou\v cek}
\address[Sou\v cek]{Matematick\'y \'ustav UK, Sokolovsk\'a 83, 186 75 Praha 8, Czech Republic}
\email{soucek@karlin.mff.cuni.cz}
\author{V\'it Tu\v cek}
\address[Tu\v cek]{Department of Mathematics, Faculty of Science, University of Zagreb, Bijeni\v cka 30, 10000 Zagreb, Croatia}
\email{}
\date{}
\thanks{P.~Pand\v zi\'c and A.~Prli\'c were supported by the QuantiXLie  Center of Excellence, a project 
cofinanced by the Croatian Government and European Union through the European Regional Development Fund - the Competitiveness and Cohesion Operational Programme 
(KK.01.1.1.01.0004). V.~Sou\v cek was supported by the grants GX19-28628X and GA24-10887S of GA\v CR. V.~Tu\v cek was supported by the QuantiXLie Center of Excellence and by the grant GX19-28628 of GA\v CR.}
\subjclass[2010]{primary: 22E47}
\keywords{}
\begin{abstract} 
{ In the 1980s, Enright, Howe and Wallach \cite{EHW} and independently Jakobsen  \cite{J} gave a complete classification of the unitary highest weight modules. In this paper we give a more direct and elementary proof of the same result for the (universal covers of the) Lie groups $Sp(2n, \mathbb{R}), SO^{*}(2n)$ and $SU(p, q)$. We also show how to describe the set of unitary highest weight modules with a given infinitesimal character.}
\end{abstract}

\maketitle

\section{Introduction}

{Let $G$ be a connected simply connected noncompact simple Lie group with Cartan involution $\Theta$ and let $K$ be the subgroup of elements fixed by $\Theta$ (then $K$ modulo center is a maximal compact subgroup of $G$ modulo center). 
Let $\frg_0=\frk_0\oplus\frp_0$ be the Cartan decomposition of the Lie algebra of $G$ corresponding to $\Theta$. As usual, we drop the subscript 0 to denote the complexifications. We assume that $(G,K)$ is a Hermitian symmetric pair, because Harish-Chandra proved n \cite {HC1} and \cite{HC2} that nontrivial unitary highest weight $(\frg,K)$-modules exist precisely when $(G,K)$ is a Hermitian symmetric pair. In case $(G,K)$ is Hermitian, the $K$-module $\frp$ decomposes into two irreducible submodules, $\frp=\frp^+\oplus\frp^-$, and each of $\frp^\pm$ is an abelian subalgebra of $\frg$. Furthermore, $\frg$ and $\frk$ have equal rank, and we choose a Cartan subalgebra $\frt$ of $\frg$ contained in $\frk$. We fix compatible choices of positive roots for $(\frg,\frt)$ and $(\frk,\frt)$, $\Delta_\frg^+\supset\Delta_\frk^+$,  such that $\frp^+$ is contained in the Borel subalgebra defined by $\Delta_\frg^+$. 
We also assume that exactly one of the simple roots is noncompact; for example, we can use the standard choices from \cite{Kn}. Our assumptions are satisfied precisely for the 
following Lie algebras $\frg_0$:  $\mathfrak{sp}(2n, \mathbb{R})$, $\mathfrak{so}^{*}(2n)$, 	$\mathfrak{su}(p, q)$, $\mathfrak{so}(2, n)$,  E III (also called $\mathfrak{e}_{6(-14)}$), and E VII (also called $\mathfrak{e}_{7(-25)}$). In this paper we will discuss the cases $\frg_0 =\mathfrak{sp}(2n, \mathbb{R})$, $\frg_0= \mathfrak{so}^{*}(2n)$, and $\frg_0=\mathfrak{su}(p, q)$. The rest of the cases will be treated in our forthcoming paper \cite{PPSST}.

A representation $M$ of $\frg$ is called a highest weight module if it is generated by a weight vector that is annihilated by the action of all positive root spaces in $\frg$. If $M$ is in addition a $(\frg,K)$-module, then the highest weight must be the highest weight of a $K$-type, hence it is dominant and integral for $\frk$. It is well known that
for each such weight $\lambda$, there is a unique irreducible $M$ as above with highest weight $\lambda$. $M$ can be constructed as the unique irreducible quotient of the generalized Verma modules with highest weight $\la$.

For $\lambda \in \frt^{*}$, the generalized Verma module $N(\lambda)$ is
\eq
\label{gen verma}
N(\lambda)=U(\frg)\otimes_{U(\frk\oplus\frp^+)} F_\la\cong F_{\lambda} \otimes S(\mathfrak{p}^{-}),
\eeq
where  $F_{\lambda}$ is the irreducible finite-dimensional $\mathfrak{k}$-module with highest weight $\lambda$. We denote by $L(\lambda)$ the irreducible quotient of $N(\lambda)$. 

Let us recall that a $(\frg, K)$-module is called unitarizable or unitary if it is isomorphic to the $(\mathfrak{g},K)$-module of $K$-finite vectors in a unitary representation of $G$. We would like to determine for which $\lambda \in \frt^{*}$ the module $L(\lambda)$ is unitary.
We already said that $\lambda$ must be $\Delta_{\mathfrak{k}}^{+}$-dominant integral. If we want $L(\lambda)$ to be unitary, then $\lambda$ must also be a real weight (i.e., contained in the real span of roots). From now on we will assume that $\lambda \in \frt^{*}$ is $\Delta_{\mathfrak{k}}^{+}$-dominant integral and real. 

A classification of unitary highest weight modules was started by  \cite{P2} and achieved in \cite{EHW} and independently in \cite{J}. 
 Later it was realized in \cite{DES} that the maximal submodule of $N(\lambda)$ is actually generated by one highest weight vector whose weight is related to a so called PRV component. This fact was used in \cite{EJ} together with a tensoring trick by \cite{J} for an alternative and uniform proof of the classification. We return to the original approach of \cite{P2} and use an extended tensoring trick to obtain uniform classification by quite elementary means. 
 
 To explain the approach of \cite{EHW} and the difference of our approach, we first recall the so called Dirac inequality. This inequality is essential in \cite{EHW}, but used in a minimal fashion. On the other hand, we use it to full extent.

Dirac operators were introduced into representation theory of real reductive groups by R.~Parthasarathy in the 1970s \cite{P1}. He successfully used them to construct most of the discrete series representations as sections of certain spinor bundles on the homogeneous space $G/K$. (His $G$ and $K$ are more general than ours: $G$ is semisimple connected with finite center, $K$ is a maximal compact subgroup, and $G$ and $K$ have equal rank.)

A byproduct of Parthasarathy's analysis was the above mentioned Dirac inequality \cite{P2}, which turned out to be a very useful necessary condition for unitarity of representations, and got to be used in several partial classifications of unitary representations.

To define Partahasarathy's Dirac operator, let $B$ be the Killing form on $\frg$; it is a nondegenerate invariant symmetric bilinear form. Let $C(\frp)$ be the Clifford algebra of $\frp$ with respect to $B$: it is the associative algebra with 1, generated by $\frp$, with relations 
\[
  xy+yx=-2B(x,y),\qquad x,y\in\frp.
\]
Let $b_i$ be any basis of $\frp$ and let $d_i$ be the dual basis with respect to $B$. The Dirac operator attached to the pair $(\frg,\frk)$ is 
\[
  D=\sum_i b_i\otimes d_i\qquad \in U(\frg)\otimes C(\frp).
\]
It is easy to see that $D$ is independent of the choice of basis $b_i$ and $K$-invariant. Moreover, $D^2$ is the spin Laplacean (Parthasarathy \cite{P1}):
\eq
\label{D squared}  
D^2=-(\Cas_\frg\otimes 1+\|\rho\|^2)+(\Cas_{\frk_\Delta}+\|\rho_\frk\|^2).
\eeq
Here $\Cas_\frg$, $\Cas_{\frk_\Delta}$ are the Casimir elements of $U(\frg)$, $U(\frk_\Delta)$, where 
$\frk_\Delta$ is the diagonal copy of $\frk$ in $U(\frg)\otimes C(\frp)$ defined by 
\[
  \frk\hookrightarrow\frg\hookrightarrow U(\frg)\quad\text{and}\quad\frk\to\frso(\frp)\hookrightarrow C(\frp). 
\]
$\rho$ and $\rho_\frk$ are as usual the half sums of positive roots for $\frg$ respectively $\frk$.

It will be convenient to choose the bases $b_i$ and $d_i$ in the following special way. Let $e_1,\dots,e_r$ be a basis of $\frp^+$ consisting of root vectors and let $f_1,\dots,f_r$ be the basis of $\frp^-$ such that $B(e_i,f_j)=\delta_{ij}$. Then $b_1,\dots,b_{2r}=e_1,\dots,e_r,f_1,\dots,f_r$ is a basis of $\frp$ with dual basis $d_1,\dots,d_{2r}=f_1,\dots,f_r,e_1,\dots,e_r$. Hence we can write
\eq
\label{D ei}
D=e_1\otimes f_1+\dots +e_r\otimes f_r+f_1\otimes e_1+\dots + f_r\otimes e_r.
\eeq

If $M$ is a $(\frg,K)$-module, then $D$ and $D^2$ act on $M\otimes S$, where $S$ is the spin module for $C(\frp)$. $S$ can be  constructed as $S=\bigwedge\frp^-$ with $\frp^-\subset C(\frp)$ acting by wedging, and $\frp^+\subset C(\frp)$ acting by contractions.

{ An inner product on the spin module $S$ is defined on $\frp^-$ by 
$$
\langle s, s'\rangle = B(s, \bar{s'}), \quad s, s' \in \frp^-,
$$
where the bar denotes the conjugation of $\frp$ with respect to $\frp_0$, and extended to all of $S$ in the usual way, using the determinant. See \cite[2.3.9]{HP2} for more details.
By \cite[Proposition 2.3.10]{HP2}, the elements of $\frp_0$ act on $S$ by skew-symmetric operators.

If $M$ is a unitary $(\frg,K)$-module, then it has an inner product such that the elements of $\frg_0$ (and hence of $\frp_0$) act on $M$ by skew Hermitian operators. We can now tensor the inner product on $M$ with the above defined inner product on $S$ to obtain an inner product on $M\otimes S$, with the property that $D$ is self adjoint with respect to this inner product. It follows that $D^2\geq 0$; this is Parthasarathy's Dirac inequality mentioned above.


Dirac inequality can be rewritten in more concrete terms using the formula \eqref{D squared} for $D^2$ and the relation between Casimir actions and infinitesimal characters. If $\Lambda$ is the infinitesimal character of $M$, and if $\tau$ is the highest weight of a $K$-type appearing in $M \otimes S$, then $D^2$ acts on the $\tau$-isotypic component of $M\otimes S$  by the scalar
\eq
\label{D squared bis}
\| \tau + \rho_{\frk}\|^2 - \| \Lambda\|^2.
\eeq
The Dirac inequality is thus equivalent to
\eq
\label{dir ind on tau}
\| \tau + \rho_{\frk}\|^2\geq \| \Lambda\|^2,
\eeq
 for any $\tau$ as above.

 In the cases we are interested in, $M=L(\la)$ or $M=N(\la)$, the infinitesimal character is $\La=\la+\rho$. We can take any $K$-type $\mu$ of $L(\la)$ and tensor it with the one-dimensional $K$-type $\bbC 1$ of the spin module $S$. The weight of $1\in S$ is $\rho_n$, the sum of noncompact positive roots (i.e., the roots of $\frp^+$), so it follows that $F_\mu\otimes 1$ is irreducible with highest weight $\mu+\rho_n$. 
 On this $K$-type of $L(\la)\otimes S$ our inequality \eqref{dir ind on tau} becomes
 \eq
\label{ehw di}
\|\mu+\rho\|\geq \|\la+\rho\|,
\eeq
and we see that if $L(\la)$ is unitary, then \eqref{ehw di} must hold for any $K$-type $\mu$ of $L(\la)$. 

This inequality is at the heart of \cite{EHW}. In fact, they prove \cite[Proposition 3.9]{EHW} that $L(\la)$ is unitary if and only if \eqref{ehw di} holds strictly for any $K$-type $\mu\neq\la$ of $L(\la)$. This criterion, while being relatively easy to prove, is difficult to use because it is hard to understand the $K$-type structure of $L(\la)$. (Even the $K$-type structure of $N(\la)$ is hard to understand explicitly because it involves tensoring, and passing to the quotient is an additional difficulty.) As a consequence, \cite{EHW} are led to prove a number of complicated statements about roots and weights related to the issue of $K$-types.

We do use \cite[Proposition 3.9]{EHW}, but only to prove unitarity of $N(\la)$ in the continuous part of the classification (for this, \cite{EHW} have another argument, which is also satisfactory; see below). For the most part however, we use instead only the lowest $K$-type $F_\la$ of $L(\la)$, but we tensor it with other $K$-types of the spin module $S$, most notably with $\frp^-\subset S=\bigwedge\frp^-$. In our approach, the basic Dirac inequality is obtained using the Parthasarthy-Ranga Rao-Varadarajan (PRV) component of $F_{\la}\otimes\frp^-$. 
The PRV component of a tensor product is described as follows.

Let $W_\frk$ denote the Weyl group of $(\frk,\frt)$. For any $\nu$ in $\frt^*$, we denote by $\nu^+$ the unique dominant $W_\frk$-conjugate of $\la$. 

\begin{prop}[\cite{PRV}]
\label{prv} 
With the above notation, let $F_\mu,F_\nu$ be finite-dimensional $\frk$-modules with highest weights $\mu,\nu$.  let $\nu^-$ be the lowest weight of $F_\nu$, and let $\tau=(\mu+\nu^-)^+$. Then $F_\tau$ appears in $F_\mu\otimes F_\nu$, with multiplicity one. Moreover, for any $F_\sigma$ appearing in $F_\mu\otimes F_\nu$,
\eq
\label{ineq prv}
\|\sigma+\rho_\frk\|^2\geq \|\tau+\rho_\frk\|^2,
\eeq
with equality attained if and only if $\sigma=\tau$.
\end{prop}

Since the lowest weight of $\frp^-$ is $-\beta$ where $\beta$ is the highest noncompact root, the Dirac inequality \eqref{dir ind on tau} for the PRV component of $F_\la\otimes\frp^-\subset L(\la)\otimes S$ is
\[
\|(\la-\beta)^+ +\rho_\frk\|^2\geq\|\la+\rho\|^2,
\]
and this is the basic inequality we are working with in this paper. In each of the cases we describe the root $\gamma$ such that 
\[
(\la-\beta)^+=\la-\gamma;
\]
the exact value of $\gamma$ depends on certain properties of $\la$. It turns out that the information about $\gamma$ is equivalent to the information about the root systems $Q$ and $R$ used by \cite{EHW} to describe the classification.

Importance of the PRV components is explained by the following corollary that was first noted in \cite{HPP}.

\begin{cor}
\label{min prv}
    Let $\mu$ be a $K$-type of $L(\la)$ and let $\nu$ be a $K$-type of the spin module $S$. Suppose that the Dirac inequality \eqref{dir ind on tau} holds for $\tau$ equal to the PRV component of $F_\mu\otimes F_\nu\subset L(\la)\otimes S$. Then \eqref{dir ind on tau} holds for any $K$-type $\tau$ contained in $F_\mu\otimes F_\nu$. 
\end{cor}
\pf This is clear from the characterization of the PRV component as the component $\tau$ with minimal $\|\tau+\rho_\frk\|$.
\epf

Corollary \ref{min prv} enables us to consider only the PRV component of each tensor product, which is mush easier then considering all components of tensor products.

The only candidate for the inner product on $L(\la)$ is the so called Shapovalov form corresponding to the conjugate involutive antiautomorphism $\sigma$ of $U(\frg)$ such that the restriction of $\sigma$ to the real form $\frg_0$ is $-\Id$. The Shapovalov form is a sesquilinear form with appropriate symmetry property
\[
\langle x v, u \rangle = \langle v, \sigma(x) u \rangle
\]
that exists for any highest weight module. It has an explicit form on $N(\la)$ and descends to a nondegenerate form on $L(\la).$ The module $L(\la)$ is unitary if and only if this induced form is positive definite.

To define the Shapovalov form, let $\frn^+$ respectively $\frn^-$ be the sum of the root spaces for positive respectively negative roots in $\Delta_\frg$. Then, thanks to our choices of Cartan subalgebra and positive roots, the restriction of $\sigma$ to $\frg$  sends positive root vectors to negative root vectors and is identity on $\frt$. 
$$
U(\frg) = U(\frt) \oplus (\frn^- U(\frg) + U(\frg) \frn^{+}).
$$
Let $P$ be the projection $U(\frg) \to U( \frt)$ along $\frn^- U(\frg) + U(\frg) \frn^{+}$.  With this we can define the universal Shapovalov form by
\[
\langle x, y \rangle = P(\sigma(x)y).
\]

For $\lambda \in \frt^*$, the Shapovalov form on the Verma module $U(\frg) \otimes_{U(\frt \oplus \frn^{+})} \mathbb{C}_{\lambda}$, with a highest weight vector $v_{\lambda}$, is defined by
$$
\langle x v_{\la}, y v_{\la} \rangle = P(\sigma (x) y)(\lambda).
$$
Here we identify $U(\frt)=S(\frt)$ with polynomials on $\frt^*$ to evaluate the polynomial $P(\sigma (x) y)$ on $\lambda$. 
The maximal submodule of the Verma module $U(\frg) \otimes_{U(\frt \oplus \frn^{+})} \mathbb{C}_{\lambda}$ is the radical of the Shapovalov form. This Shapovalov form on the Verma module induces forms on $N(\lambda)$ and $L(\lambda)$ which we also call Shapovalov forms. 

The shape of the \cite{EHW} classification is given by certain lines in $\frt^*$ of the form 
\eq
\label{lines}
\la=\lambda_0 + z \zeta, \qquad z \in \mathbb{R}.
\eeq 
Here $\zeta\in\frt^*$ is orthogonal to $\Delta_{\frk}$ and normalized so that  
$\frac{2 \langle \zeta, \beta \rangle}{\langle \beta, \beta \rangle } = 1$, 
where $\beta$ as before denotes the unique maximal noncompact root of $\Delta_\frg^+$. For a fixed $\la$, $\la_0$ is defined as the point on the line $\la+ z\zeta$, $z\in\bbR$,  such that $\langle \lambda_0 + \rho , \beta \rangle = 0$. The line is then renormalized so that $\la$ varies as in \eqref{lines}.

It is known from the work of Harish-Chandra that if $\la$ is on the line \eqref{lines}, for sufficiently negative $z$ the module $N(\la)$ is the $(\frg,K)$-module of a holomorphic disrete series representation. In particular it is irreducible and unitary, and its Shapovalov form is thus positive definite. 

Let $a$ be the smallest real number such that $N(\lambda_0 + a \zeta)$ is reducible. The corresponding $\la$ is called the first reduction point. It follows by continuity that the Shapovalov forms on $N(\lambda_0 + z \zeta)$, $z<a$, are all positive definite and hence these modules are unitary. 
Therefore the Shapovalov form on $N(\lambda_0 + a \zeta)$ is positive semidefinite  and it induces a positive definite form on $L(\lambda_0 + a \zeta)$. Hence $L(\lambda_0 + a \zeta)$ is also unitary. \cite{EHW} use this argument to prove unitarity of $N(\lambda_0 + z \zeta)$, $z<a$, and also of $L(\la_0+a\zeta)$. 
See \cite[Proposition 3.1]{EHW}. 
To identify the point $a$ they use Jantzen's criterion for irreducibility of generalized Verma modules. We obtain irreducibility and unitarity of $N(\lambda_0 + z \zeta)$, $z<a$, directly from Dirac inequality, using Corollary \ref{cor unit nonunit} below. It then follows that $\lambda_0 + a \zeta$ is the first reduction point; see Remark \ref{rem 1st red pt}.

The irreducible unitary modules $N(\lambda_0 + z \zeta)$, $z<a$, form the continuous part of the set of the unitary points on the line \eqref{lines}. In addition, there is a finite number of discrete unitary points on the line. One of these discrete points (possibly the only one) is obtained for $z=a$. 

\cite{EHW} construct the (unitary) representations in the discrete part of the classification by using Howe's theory of dual pairs and theta correspondences.  (In some cases they need additional considerations.)  
We construct these representations by a tensoring argument, which is similar in spirit but simpler and more elementary. We first construct some basic representations starting from the Weil representation of $\frsp(2n,\bbR)$. Then we construct other  representations in the discrete part as certain submodules of the tensor products of the basic representations. This construction is described at the end of Section \ref{sec prelim}, and then in more detail for each of the cases we are considering. 

In most cases we can use this tensoring technique to prove unitarity of $L(\la)$ at all discrete points, including the first reduction point $z=a$. In some cases we get from tensoring all the discrete points except for $z=a$; this happens in this paper for some cases when $\frg_0=\frso^*(2n)$. In those cases we use the above continuity argument for unitarity at $z=a$. We record the statement we need:
\begin{lem}
    \label{continuity}
    Assume that on the line \eqref{lines} the modules $N(\la_0+z\zeta)$ are irreducible and unitary for $z<a$. Then $L(\la_0+a\zeta)$ is unitary. \qed
\end{lem}


To talk about $K$-types of $N(\la)$ and $L(\la)$, one first needs to describe the $K$-types of $S(\frp^-)$. These $K$-types 
are well known from the work of Schmid 
\cite{S} and we call them the Schmid modules. It is customary to describe the $K$-types of $S(\frp^+)$ and then to pass to $S(\frp^-)$ by duality; if $F_s$ is a $K$-type of $S(\frp^+)$ with highest weight $s$, then the corresponding $K$-type of $S(\frp^-)$ is $F_s^*=F_{-s}$, with lowest weight $-s$. We identify the modules $F_{-s}$ with the weights $s$ and talk about Schmid modules $s$.

The basic Schmid modules $s_1,\dots,s_r$ are obtained from the Harish-Chandra system of strongly orthogonal positive roots; these are described for example in \cite[1.4]{EJ}, and we give them explicitly for each of the cases we consider.
If the strongly orthogonal roots are $\gamma_1,\dots,\gamma_r$, then the basic Schmid modules are
\eq
\label{basic schmid}
s_i=\gamma_1+\ldots+\gamma_i,\qquad i=1,\dots,r.
\eeq
The first basic Schmid module is $s_1=\gamma_1=\beta$, the highest noncompact root; the corresponding $K$-type in $S(\frp^-)$ is $\frp^-$. 

A general Schmid module $s$ is a nonnegative integer combination of the basic Schmid modules:
\eq
\label{general schmid}
s=a_1 s_1 + \ldots + a_r s_r,\qquad a_1,\dots,a_r\in \bbZ_+.
\eeq
The $K$-module $S(\frp^-)$ decomposes as the direct sum of all Schmid modules, with multiplicities equal to 1.

It is useful to introduce the concept of level: a $K$-type of $S(\frp^-)$ has level $n$
if it appears in $S^n(\frp^-)$. The level of the basic Schmid module $s_i$ is $i$, and the level of a general Schmid module \eqref{general schmid} is $\sum ia_i$. 

The level of a $K$-type $\mu$ of $N(\lambda)$ is the number $n$ such that $\mu$ is contained in $F_{\lambda} \otimes S^n(\frp^-)$. 
(We note that \cite{EJ} use the word level for a different notion.)

We use the concept of level in Corollary \ref{cor unit nonunit} where we give a criterion for non-unitarity of $L(\la)$ in certain cases. This is later used in each of the cases to obtain all non-unitary points on each of the lines \eqref{lines}.

\bigskip 

\section{Some preliminary facts}
\label{sec prelim}

Some of the more involved computations needed to prove the relevant Dirac inequalities in \cite{PPST1} and \cite{PPST2} use the so called generalized PRV components of the tensor product. The following proposition was conjectured  by Parthasarthy-Ranga Rao-Varadarajan \cite{PRV} and proved by Kumar \cite{Ku} and independently by Mathieu \cite{M}.

\begin{prop}
\label{prop gen prv}
Let $F_\mu$ and $F_\nu$ be the irreducible finite-dimensional $K$-modules with highest weights $\mu$ respectively $\nu$. Let $w\in W_\frk$ and let $\tau=(\mu+w\nu)^+$ be the $\frk$-dominant conjugate of $\mu+w\nu$. Then the $K$-module $F_\tau$ with highest weight $\tau$ appears in $F_\mu\otimes F_\nu$ (with multiplicity at least one).\qed
\end{prop}

We note that Proposition \ref{prv} and Proposition \ref{prop gen prv} are originally stated and proved for modules over a semisimple Lie algebra, however the results easily extend to modules over a reductive Lie algebra $\frk$. In particular, the inequality \eqref{ineq prv} remains unchanged. Namely, in the reductive case we can write $\mu=\mu'+\eta$, $\nu=\nu'+\chi$, where $\mu',\nu'$ are the highest weights of $F_\mu,F_\nu$ restricted to the semisimple part of $\frk$ and $\eta,\chi$ are the central characters. The central characters are orthogonal to the roots of $\frk$, while the weights $\mu',\nu'$ are in the span of the roots of $\frk$. Then all the components of $F_\mu\otimes F_\nu$ have central character $\theta=\eta+\chi$, and we can write $\sigma=\sigma'+\theta$, 
$\tau=\tau'+\theta$, with $\theta$ orthogonal to both $\sigma'$ and $\tau'$. It is now clear that $\|\sigma'+\rho_\frk\|^2\geq \|\tau'+\rho_\frk\|^2$ is equivalent to
\[
\|\sigma+\rho_\frk\|^2=\|\sigma'+\rho_\frk\|^2+\|\theta\|^2\geq \|\tau'+\rho_\frk\|^2+\|\theta\|^2=\|\tau+\rho_\frk\|^2.
\]
By the same argument, if $(\frg,\frk)$ is a Hermitian symmetric pair, and if $F_\mu,F_\nu$ are finite-dimensional $\frk$-modules as above, the inequality \eqref{ineq prv} is equivalent to
\eq
\label{ineq prv2}
\|\sigma+\rho\|^2\geq \|\tau+\rho\|^2.
\eeq
Namely, $\rho_n=\rho-\rho_\frk$ is orthogonal to the roots of $\frk$.

\begin{cor}
\label{gen prv}
Let $(\frg,\frk)$ be a Hermitian symmetric pair. Let $\mu$ and $\nu$ be  dominant integral $\frk$-weights. Let $w_1,w_2\in W_\frk$. Then
\[
\|(w_1\mu-w_2\nu)^++\rho\|^2\geq \|(\mu-\nu)^++\rho\|^2.
\]
\end{cor}
\pf
Note that $(\mu-\nu)^+$ is the highest weight of the PRV component of $F_\mu\otimes F_\nu^*$. On the other hand, $(w_1\mu-w_2\nu)^+=(\mu-w_1^{-1}w_2\nu)^+$ is the highest weight of a generalized PRV component of $F_\mu\otimes F_\nu^*$, corresponding to the extremal weight $-w_1^{-1}w_2\nu$ of $F_\nu^*$. By Proposition \ref{prop gen prv}, this generalized PRV component appears in the tensor product. The required inequality now follows from \eqref{ineq prv2}.
\epf

\begin{lem} 
\label{dir eq max sub}
Let $M(\la)$ be the maximal submodule of the generalized Verma module $N(\la)$, and let $F_\mu\subset M(\la)$ be a maximal 
$K$-type of $M(\la)$, with highest weight $\mu$. Then
\[
D^2=0\qquad\text{on}\quad F_\mu\otimes 1\subset N(\la)\otimes S.
\]
Equivalently, 
\[
\|\mu+\rho\|^2=\|\la+\rho\|^2.
\]
\end{lem}
\pf
The two statements are equivalent because $N(\la)$ has infinitesimal character $\la+\rho$ and $F_\mu\otimes 1$ has highest weight $\mu+\rho_n$, so the action of $D^2$ on $F_\mu\otimes 1$, given by the scalar \eqref{D squared bis}, is equal to
\[
\|\mu+\rho\|^2-\|\la+\rho\|^2.
\]
(See also the discussion above \eqref{ehw di}.)

On the other hand, if $D$ is written as in \eqref{D ei}, then 
\[
D(F_\mu\otimes 1) = \sum (e_iF_\mu\otimes f_i\cdot 1 + f_iF_\mu\otimes e_i\cdot 1) =0,
\]
since $e_i F_\mu=e_i\cdot 1=0$.

Alternatively, one can note that $F_\mu$ generates a 
$\frg$-submodule of $N(\la)$ of highest weight $\mu$, which must have the same infinitesimal character as $N(\la)$, so $\mu+\rho$ and $\la+\rho$ are conjugate by the Weyl group of $\frg$, hence have the same norm.
\epf


The following results were used in \cite{PPST1} and \cite{PPST2}, but their proofs fit better into this paper.

Note that since $N(\la)=F_\la\otimes S(\frp^-)$, the $K$-type with highest weight $(\la-s)^+$ appears in $N(\la)$ for any Schmid module $s$. Namely, this $K$-type is the PRV component of $F_\la\otimes F_{-s}$.

\begin{cor}
\label{cor unit nonunit}
(1) Let $s_0$ be a Schmid module such that 
the strict Dirac inequality \eqref{ehw di} for $\mu=(\la-s)^+$, i.e.,  
\eq
\label{strict di s}
\|(\la-s)^++\rho\|^2> \|\la+\rho\|^2
\eeq
holds for any Schmid module $s$ of strictly lower level than $s_0$, and such that
\[
\|(\la-s_0)^++\rho\|^2< \|\la+\rho\|^2.
\]
Then $L(\la)$ is not unitary.

(2) If \eqref{strict di s} 
holds for all Schmid modules $s$, then $N(\la)$ is irreducible and unitary.
\end{cor}
\pf
(1) By assumption, the Dirac inequality fails for the $K$-type 
$(\lambda-s_0)^+$ of $N(\la)$, so we only need to show that $(\lambda-s_0)^+$ is a $K$-type of $L(\lambda)$ and not just of $N(\lambda)$. 

Suppose $(\la-s_0)^+$ is a $K$-type of the maximal submodule $M(\la)$ of $N(\la)$. Then, by Lemma \ref{dir eq max sub}, it can not be a maximal $K$-type of $M(\la)$ and so the $U(\frg)$-action generates from it a $K$-type of $M(\la)$ of lower level. It follows that there must exist a maximal $K$-type $\nu$ of $M(\la)$ with lower level than that of $(\la-s_0)^+$. On this maximal $K$-type we have $\|\nu+\rho\|^2 - \|\la+\rho\|^2 = 0$ by Lemma \ref{dir eq max sub}. At the same time, due to Proposition \ref{prv}, there exists a Schmid module $s$ such that $\|\nu+\rho\|^2 - \|\la+\rho\|^2$ is bounded below by the value of the corresponding expression for $(\la - s)^+$ which is positive by \eqref{strict di s}.

(2)
The assumption implies that 
\eq
\label{strict di verma}
\|\mu+\rho\|^2>\|\la+\rho\|^2
\eeq
for any $K$-type $\mu\neq\la$ of the generalized Verma module $N(\la)$. By Lemma \ref{dir eq max sub}, it follows that $N(\la)$ is irreducible. 

Unitarity of $N(\la)$ follows from the strict Dirac inequality \eqref{strict di verma} by \cite[Proposition 3.9]{EHW}.
\epf

\begin{rem}
\label{rem 1st red pt}
{\rm
In each of the cases, we will use Corollary \ref{cor unit nonunit} to see that modules $N(\la_0+z\zeta)$ are irreducible and unitary for $z<a$ for a certain value of $a$, and that $L(\la_0+z\zeta)$ is not unitary for $z>a$ close to $a$ (more precisely, for $z$ between $a$ and the next discrete point in the classification).

This implies that $N(\la_0+a\zeta)$ can not be irreducible, so $\la_0+a\zeta$ is the first reduction point on the line $\la_0+z\zeta$. Namely the Shapovalov form on $N(\la_0+z\zeta)$ is positive definite for $z<a$ and has a strictly negative eigenvalue for $z>a$ close to $a$. Thus by continuity it must have a nontrivial radical at $z=a$, and this radical is the maximal submodule of $N(\la_0+a\zeta)$.
}
\end{rem}

 We will construct the most interesting unitary highest weight modules, those in the discrete part of the classification, by tensoring some basic modules, starting from the Weil representation. For this construction, the following proposition will be crucial.
 
\begin{prop}\label{dot}
    Let $V_1$ and $V_2$ be irreducible unitary highest weight $(\frg,K)$-modules 	with top $K$-types $F_1$ and $F_2.$ Then:

    (1) $V_1\otimes V_2$ is a unitary $(\frg,K)$-module;

    (2) Let $F$ be an irreducible $K$-submodule of $V_1\otimes V_2$ annihilated by $\frp^+$.  Then $U(\frp^-)F$ is an irreducible unitary highest weight module;

    (3) 	Let $F$ be an irreducible $K$-submodule  of  $F_1\otimes F_2$.  Then $U(\frp^-)F$ is an irreducible unitary highest weight module.
\end{prop}

\pf
(1) Note that any irreducible highest weight $(\frg,K)$-module contains only finitely many $K$-types with a given central character. This implies that $V_1\otimes V_2$ is a $(\frg,K)$-module. Unitarity of $V_1\otimes V_2$ is clear, as we can take the tensor product of the inner products on $V_1$ and $V_2$.

(2) $U(\frp^-)F$ is unitary since it is a submodule of the module $V_1\otimes V_2$ which is unitary by (1). Suppose $U(\frp^-)F$ is not irreducible. Then it contains a proper nonzero submodule $M$, which must have a highest weight vector outside of $F$. By unitarity, $U(\frp^-)F=M\oplus M^\perp$, where $M^\perp$ is a submodule containing $F$. But then $M^\perp$ contains all of $U(\frp^-)F$, a contradiction. Thus $U(\frp^-)F$ is irreducible.

(3) is a special case of (2), since $F_1\otimes F_2$ is annihilated by $\frp^+$.
\epf

The PRV components in a tensor product of $K$-modules can be used to construct a useful product operation for irreducible unitary highest weight $(\frg,K)$ modules.

\begin{df}\label{bullet}
		Let $V_1, V_2, F_1, F_2$ be as in Proposition \ref{dot}. If $F$ is  the PRV component of $F_1\otimes F_2,$ 
	we shall call the module $U(\frp^-)F$ the PRV product of $V_1$ and $V_2$ and we denote it by $V_1\bul V_2.$	
\end{df}}

Our construction can be viewed as an analogue of the construction of finite-dimensional modules of a simple Lie algebra $\frg$, where one first considers the standard module. Then one constructs the fundamental representations by taking the exterior powers of the standard module, and possibly adds the spin module(s). These exterior powers are the PRV components of the tensor products of the standard module with itself. Finally, one takes the Cartan components of the tensor products of the fundamental representations to obtain all finite-dimensional representations.

The cone structure of the discrete part of the set of all irreducible unitary highest weight
modules was described and discussed in \cite{DES}.
We shall see below that the complete cone structure for a given case can be obtained by first constructing the cones depending on one parameter  only (the basic representations), and then by taking their
PRV products. 
}

\section{Classification for $\frsp(2n,\bbR)$}

Let $\lambda$ be the highest weight of an irreducible $(\frg,K)$-module $L(\lambda)$. Recall that $L(\lambda)$ is the simple quotient of the generalized Verma module $N(\lambda)$ with highest weight $\la$. The infinitesimal character of $L(\la)$ and $N(\la)$ is represented by the parameter $\lambda+\rho$, which must be dominant regular integral for $\frk$. In other words, the coordinates of $\la+\rho$ must be strictly decreasing and their differences must be integers.

The basic inequality that $\la$ must satisfy if $L(\la)$ is unitary is
\eq
\label{basic di sp}
\|(\la-\beta)^++\rho\|^2\geq\|\la+\rho\|^2,
\eeq
where $\beta=2\eps_1$ is the highest root. 
To understand this inequality better, let $q\leq r$ be integers in $[1,n]$ such that
\eq
\label{lambda sp}
\lambda=(\underbrace{\lambda_1,\dots,\lambda_1}_q,\underbrace{\lambda_1-1,\dots,\lambda_1-1}_{r-q},\lambda_{r+1},\dots,\lambda_n),
\eeq
with $\lambda_1-2\geq \lambda_{r+1}\geq\dots\geq\lambda_n$.
Then 
\[
(\la-\beta)^+=(\underbrace{\lambda_1,\dots,\lambda_1}_{q-1},\underbrace{\lambda_1-1,\dots,\lambda_1-1}_{r-q},\la_1-2,\lambda_{r+1},\dots,\lambda_n)=\la-(\eps_q+\eps_r).
\]
The root $\gamma=\eps_q+\eps_r$, or the integers $q$ and $r$, contain the same information about $\la$ as the root systems $Q$ and $R$ of \cite{EHW}. In fact, it is easy to see that $Q$ is the root system of the subalgebra $\frsp(2q,\bbR)$ of $\frg$ built on the first $q$ coordinates, while $R$ is the root system of the subalgebra $\frsp(2r,\bbR)$ of $\frg$ built on the first $r$ coordinates.

As shown in \cite{PPST1}, the inequality \eqref{basic di sp} is equivalent to

\eq
\label{basic cond sp}
\la_1\leq -n+\frac{r+q}{2}.
\eeq

Recall that by a result of Schmid \cite{S}, the $K$-module  
$S(\frp^-)$ is multiplicity free, and its $K$-types are the irreducible $K$-modules $F_{-s}$ with lowest weights
$-s$, where 
\eq
\label{gen schmid}
s=(2b_1,\dots,2b_n),\qquad b_i\in\bbZ,\quad b_1\geq\dots\geq b_n\geq 0.
\eeq
We will refer to such $F_{-s}$, or to $s$, as Schmid modules; they are nonnegative integer combinations of the basic Schmid modules
\[
s_i=(\underbrace{2,\dots,2}_i,0,\dots,0),\qquad i=1,\dots,n.
\] 
Since $N(\lambda)=F_\lambda\otimes S(\frp^-)$, the $K$-types of $N(\lambda)$ are irreducible summands of various $F_\lambda\otimes F_{-s}$, where $F_{-s}$ is as above. We want to see for which of them the Dirac inequality holds. Note that by Proposition \ref{prv}, if for some $s$ the PRV component of $F_\lambda\otimes F_{-s}$ satisfies the Dirac inequality, i.e., if
\eq
\label{din1}
\|(\la-s)^++\rho\|^2\geq \|\la+\rho\|^2,
\eeq
then the same is true for all components of $F_\lambda\otimes F_{-s}$. Thus it will be enough for us to work with the PRV components. We start by examining what happens for the first $q$ basic Schmid modules $s=s_1,\dots,s_q$. Since $s_1=\beta$ (the highest root), we already know that the Dirac inequality for $s_1$ is satisfied if and only if \eqref{basic cond sp} holds. More generally, as was shown in \cite{PPST1}, \eqref{din1} holds for $s=s_i$, $i\in\{1,\dots,q\}$, if and only if  
\eq
\label{din2}
\lambda_1\leq -n+\frac{r+q-i+1}{2}.
\eeq
Moreover,  \eqref{din1} holds strictly for $s_i$ if and only if \eqref{din2} holds strictly.
Note that for $i=1$, \eqref{din2} is exactly our basic inequality \eqref{basic cond sp}, while for $i=q$ we get 
\[
\lambda_1\leq -n+\frac{r+1}{2}.
\]
It turns out that the critical values for $\la_1$ are exactly those appearing in \eqref{din2}, i.e., 
\[
-n+\frac{r+1}{2},-n+\frac{r+2}{2},\dots,-n+\frac{r+q}{2}.
\]
In Theorem \ref{thm sp} below we will show that if $\la_1$ is equal to one of these critical values, then $L(\la)$ is unitary. All other $\la$ are handled by the following theorem which follows from the results of \cite{PPST1}.

\begin{thm}
\label{unit nonunit sp}
Let $\la$ be as in \eqref{lambda sp}. Then:
\begin{enumerate}
\item If 
\[
\la_1>-n+\frac{r+q}{2},
\]
then $L(\la)$ is not unitary.
\item If for some integer $i\in[1,q-1]$ 
\eq
\label{gaps}
-n+\frac{r+q-i}{2}<\lambda_1< -n+\frac{r+q-i+1}{2},
\eeq
then $L(\la)$ is not unitary.
\item If
\eq
\label{cont}
\lambda_1< -n+\frac{r+1}{2},
\eeq
then $L(\lambda)=N(\lambda)$ and it is unitary. 
\end{enumerate}
\end{thm}
\pf
(1) is immediate from our basic inequality \eqref{basic cond sp}.

(2) follows from Theorem 3.1.(1)
of \cite{PPST1}, from Corollary \ref{cor unit nonunit}.(1), and from the obvious fact that a Schmid module $s$ of lower level than $s_{i +1}$ can have at most $i$ nonzero components. 

(3) follows from Theorem 3.1.(2) of \cite{PPST1} and from  Corollary \ref{cor unit nonunit}.(2). 
\epf

We say that the $\la$ satisfying \eqref{cont} are in the continuous part of its line, where the line is meant in the sense of \cite{EHW}, i.e., $\lambda,\la'$ are on the same line if $\la-\la'$ has all coordinates equal to each other.
\smallskip

To complete the classification, we prove the following theorem.

\begin{thm}\label{thm sp}
Let $1\leq \ell\leq q\leq r\leq n$ and $a_1\geq a_2\geq\ldots\geq a_{n-r}\geq 0$  
 be two sequences of integers. Let 
\eq\label{la1 sp}
\la_1=-n+\frac{q+r-\ell+1}{2}
\eeq 
and let $L_\la$ be the irreducible module with highest weight 
\eq\label{unitary sp}
\la=(\underbrace{\la_1,\ldots,\la_1}_q,\underbrace{\la_1-1,\ldots, \la_1-1}_{r-q},
\underbrace{\la_1-2-a_{n-r},\ldots,\la_1-2-a_1}_{n-r}).
\eeq
 Then $L_\la$ is unitary.
\end{thm}

The cones of \cite{DES} are easily described in terms of \eqref{unitary sp}. The vertex of the cone containing $\lambda$ is obtained by setting all $a_k$ equal to zero, and all elements of the cone are obtained by taking all possible choices for the $a_k$.

The proof of Theorem \ref{thm sp} is based  on properties of the well known unitary  Segal-Shale-Weil
(or oscillatory) representation $W$ of the Lie algebra $\frg=\frsp(2n).$  
	
Let us first recall the standard realization of  $W$ (see, e.g., \cite{H3}). The representation is defined on the vector space
$W=\bbC[z_1,\ldots,z_n]$ of polynomials in $n$ variables. Let
$\pa_i=\frac{\pa}{\pa z_i}.$

Basis elements of $\frk\subset\frg$ act as first order differential operators
\eq
-\frac{1}{2}(z_i\pa_j+\pa_jz_i);\;1\leq i,j\leq n
\eeq
The span of the set $\{-z_i\pa_i-\frac{1}{2};\;1\leq i\leq n\}$ is the Cartan subalgebra $\frt$ of both $\frk$ and $\frg$. The positive root vectors for $\frk=\fru(n)_\bbC$ are $\del_i z_j$, $1\leq i<j\leq n$; these are root vectors for the positive roots $\eps_i-\eps_j$.

Basis elements of $\frp^+\subset\frg_\bC$
acts by the second order differential operators $\pa_i\pa_j;\;1\leq i,j\leq n$; these are roots vectors for the positive roots $\eps_i+\eps_j$. 
Basis elements of $\frp^-\subset\frg_\bC$ acts by zero order differential operators
$z_iz_j;\;1\leq i,j\leq n$; these are root vectors for the negative roots $-\eps_i-\eps_j$. 

It is well known that $W$ is a unitary module, which decomposes
into a direct sum $W_+\oplus W_-$ of two irreducible highest weight submodules, which are realized
by even, respectively odd, polynomials. The highest weight of $W_+$ is 
$\la_+=(-\frac{1}{2},\ldots,-\frac{1}{2})$ and the highest weight of $W_-$ is
$\la_-=(-\frac{1}{2},\ldots,-\frac{1}{2},-\frac{3}{2}).$

\begin{df}
	Let $L(\la);\;\la\in\frt^*$ 
 denote the irreducible highest weight $(\frg,K)$-
	module with highest weight $\la.$

A basic sequence of irreducible highest weight $(\frg,K)$-modules $L_{-k},$ is defined by
\eq\label{basic}
L_{-k}:=L(\la_k);\;\la_k=(-1,\ldots,-1,-k);\;k\geq 2.
\eeq
\end{df}

We will call the modules $L_{-k}$  basic representations. As we shall see, together with the two Weil representations they can be used as building blocks for unitary highest weight modules. The more complicated modules are obtained from the basic ones and the two Weil representations by tensoring.

\begin{lem}  
The module $L_{-k}$, $k\geq 2,$ is isomorphic to an irreducible submodule of $W\otimes W_+,$
hence it is  an irreducible unitary highest weight module.
\end{lem}

\pf 
We claim that for every $m\geq 1,$ the vector 
\eq\label{sing}
w=\sum_{i=0}^{[m/2]} (-1)^i\binom{m}{2i} z_n^{m-2i}\otimes z_n^{2i}
\eeq
is a singular vector for $\frg$, i.e., it is killed by all positive $\frk$-root vectors $\del_i z_j$, $1\leq i<j\leq n$, and by all $\del_i\del_j\in\frp^+$. This will prove the lemma, since it will follow from Proposition  \ref{dot} that $w$ generates an irreducible unitary submodule of $W\otimes W_+$, and since the weight of $w$ is
\[
(-\half,\dots,-\half,-\half-(m-2i))+(-\half,\dots,-\half,-\half-2i)=(-1,\dots,-1,-1-m),
\]
the module generated by $w$ is isomorphic to $L_{-m-1}$.

It is clear that $w$ is killed by all $\del_i z_j$, $i<j$, and by all $\del_i\del_j$ except possibly by $\del_n^2$.
Furthermore, 
\[
\del_n^2 w =\sum_{j=0}^{[m/2]-1} a_jz_n^{m-2-2j}\otimes z_n^{2j},\qquad \text{for some } a_j\in\bbC.
\]
The summand $a_jz_n^{m-2-2j}\otimes z_n^{2j}$ is obtained when $\del_n^2\otimes 1$ acts on $z_n^{m-2j}\otimes z_n^{2j}$ and when $1\otimes\del_n^2$ acts on $z_n^{m-2-2j}\otimes z_n^{2j+2}$. It follows that
\[
a_j=(-1)^j\binom{m}{2j}(m-2j)(m-2j-1) + (-1)^{j+1}\binom{m}{2j+2}(2j+2)(2j+1)=0.
\]
So $\del_n^2 w=0$ and the lemma follows.
\epf

\begin{ex}
\label{ ex prv} 
{\rm
Let us show how to obtain $L(\la)=L(-4,-4,-5,-5,-7,-8)$ 
from basic representations and Weil representations using the PRV product defined in Definition \ref{bullet}. This $\la$ appears in  Theorem \ref{thm sp} for $n=6$, $q=2$, $r=4$, $\ell=3$, $a_1=2$ and $a_2=1$. Since PRV products of unitary modules are unitary, it will follow that 
$L(\la)$
is unitary.

We first note that $L_{-4}\bul L_{-5}$ has highest weight
\[
\mu=[(-1,\dots,-1,-4)+(-5,-1,\dots,-1)]^+=(-2,-2,-2,-2,-5,-6).
\]
Furthermore, $(W_-)^2=W_-\bul W_-$ has highest weight
\[
\nu=\left[\left(-\half,\dots,-\half,-\frac{3}{2}\right)+\left(-\frac{3}{2},-\half,\dots,-\half\right)\right]^+=(-1,-1,-1,-1,-2,-2).
\]
We now consider $L_\mu\bul L_\nu$; its highest weight is
\[
\eta=[(-2,\dots,-2,-5,-6)+(-2,-2,-1,\dots,-1)]^+=(-3,-3,-4,-4,-6,-7).
\]
Finally, we make the PRV product of $L_\eta$ and 
\[
(W_+)^2=W_+\bul W_+=L_{(-1,\dots,-1)}
\]
to obtain $L_\la$.
}
\end{ex}

{\it Proof of Theorem \ref{thm sp}.} 
Let
\[
L_{a_1,\ldots, a_{n-r}}=
 L_{-3-a_{n-r}}\bul (L_{-3-a_{n-r-1}}\bul(\ldots \bul(L_{-3-a_2}\bul L_{-3-a_1})\ldots)).
\]
We also define
\[
(W_-)^{r-q}=
 \underbrace{W_-\bul (W_-\bul(\ldots \bul(W_-\bul W_-)\ldots))}_{r-q}\qquad\text{and}
 \]
 \[ 
 (W_+)^{\ell-1}=
 \underbrace{W_+\bul (W_+\bul(\ldots \bul(W_+\bul W_+)\ldots))}_{\ell-1}.
\]
We are going to show that
\eq\label{decomp}
  L(\la)\cong \left( L_{a_1,\ldots a_{n-r}}\bul (W_-)^{r-q}\right)\bul (W_+)^{\ell-1}.
\eeq
This will imply the theorem since PRV products of unitary modules are unitary. 

\begin{lem}\label{W_-}
Let $ s$ be an integer with $1\leq s\leq n$.
Then the irreducible unitary highest weight module $(W_-)^s$ has
highest weight
\[
\la_s=( \underbrace{-\frac{s}{2},\ldots,-\frac{s}{2}}_{n-s},
-\frac{s}{2}-1,\ldots,-\frac{s}{2}-1 ).
\]
 \end{lem}
 
\pf
We use induction on $s$. The case $s=1$  is trivial.
 
Suppose that the statement holds for $s$ and consider
the tensor product 
$L(\la_s)\otimes W_-.$ 
The lowest weight of the top K-type
of $W_-$ is equal to $\mu=(-\frac{3}{2},-\frac{1}{2},\ldots,-\frac{1}{2}),$ 
hence the PRV component of the tensor product of top $K$-types of both factors
has the highest weight equal to 
 $$
 (\la_s+\mu)^+  = 
\underbrace{ \left(\frac{-s-3}{2},\frac{-s-1}{2}, \ldots,\frac{-s-1}{2}\right.}_{n-s},
\left.\frac{-s-3}{2},
 \ldots,\frac{-s-3}{2}\right)^+,
$$
which is equal to $\la_{s+1}.$
\epf
	
\begin{lem}\label{conus}
Fix $j$ with $1\leq j\leq n.$
Let $ a_1\geq\ldots\geq a_j\geq 0$ be a sequence of integers and consider the irreducible unitary highest weight module  
\[
L_{a_1,\ldots a_j}=
L(-3-a_j)\bul (L(-3-a_{j-1})\bul(\ldots \bul(L(-3-a_2)\bul L(-3-a_1))\ldots)).
\]
Then $L_{a_1,\ldots a_j}$ has highest weight 
\[
\la_{a_1,\ldots,a_j}=(\underbrace{-j,\ldots,-j}_{n-j},
-j-2-a_j,\ldots,-j-2-a_1).
\]
\end{lem}
 
 \pf
 We use induction on $j$. The case $j=1$  is trivial.
 
 Suppose that the Lemma holds for the sequence $a_1\geq\ldots\geq a_{j-1}$, i.e., that $L_{a_1,\ldots, a_{j-1}}$ has highest weight
 \[
 (\underbrace{-j+1,\ldots,-j+1}_{n-j+1},
-j-1-a_{j-1},\ldots,-j-1-a_1).
\]
Consider the  product $L_{-3-a_j}\bul L_{a_1,\ldots, a_{j-1}}.$ The lowest weight of the top K-type
 of $L_{a_1,\ldots, a_{j-1}}$ is equal to $(-j-1-a_1,\ldots,-j-1-a_{j-1},-j+1,\ldots,-j+1)$, so $L_{-3-a_j}\bul L_{a_1,\ldots, a_{j-1}}$  
 has highest weight 
 \[
 (-j-2-a_1,\ldots,-j-2-a_{j-1},
 -j,\ldots,-j,-j-2-a_j)^+
\] 
 which is equal to $\la_{a_1,\ldots, a_j}.$ 
 \epf
 
We now get back to the proof of Theorem \ref{thm sp}. From   Lemma \ref{W_-} we see that  
 the lowest weight of the top $K$-type of $(W_-)^{r-q}$ is 
\[
(\underbrace{\frac{q-r}{2} - 1,\ldots,\frac{q-r}{2}-1}_{r-q},
\frac{q-r}{2},\ldots,\frac{q-r}{2})
\]

Using this and Lemma \ref{conus} we see that 
the highest weight of $  L_{a_1,\ldots, a_{n-r}}\bul (W_-)^{r-q}$
is 
$$
(\underbrace{\alpha-1,\ldots,\alpha-1}_{r-q},
\underbrace{\alpha,\ldots, \alpha}_{q},
\underbrace{ \alpha-2-a_{n-r},
		\ldots,\alpha-2-a_1}_{n-r})^+,
$$
where $\alpha=-n+\frac{q+r}{2}.$

Since the lowest weight of the top $K$-type of $(W_+)^{\ell-1}$
is 
$$
\left(\frac{1-\ell}{2},\ldots,\frac{1-\ell}{2}\right),
$$
we see that $\left( L_{a_1,\ldots a_{n-r}}\bul (W_-)^{r-q}\right)\bul (W_+)^{\ell-1}$ has highest weight $\la$. This proves \eqref{decomp} and hence finishes the proof of Theorem \ref{thm sp}.
\qed

\section{Classification for $\frsu(p,q)$, $p\leq q$}

Let $\lambda$ be the highest weight of an irreducible $(\frg,K)$-module $L(\lambda)$. 

We use the quotient version for the standard coordinates of the complexified compact Cartan subalgebra of $\frsu(p,q)$ and its dual: setting $n=p+q$, we consider $n$-tuples of complex numbers modulo the $n$-tuples with all coordinates equal to each other. We note that using these coordinates does not cause any problems when working with the norms; the argumentation for that is the same as the discussion below Proposition \ref{prop gen prv}.

We set $n=p+q$ and write $\la$ in the standard coordinates as
\[
\lambda=(\lambda_1,\dots,\lambda_p\,|\,\lambda_{p+1},\dots,\lambda_n),
\]
where components $\la_1,\dots,\la_n$ satisfy
\[
\la_1\geq\dots\geq\la_p;\qquad \la_{p+1}\geq\dots\geq\la_n,
\]
and $\la_i-\la_j$ is an integer for $i, j \in \{ 1, \ldots, p\}$ and $i, j \in \{ p+1, \ldots, n\}$.

The basic inequality that $\la$ must satisfy if $L(\la)$ is unitary is
\eq
\label{basic di su}
\|(\la-\beta)^++\rho\|^2\geq\|\la+\rho\|^2,
\eeq
where $\beta=\eps_1-\eps_n$ is the highest (noncompact) root. 
To understand this inequality better, let $p'\leq p$ and $q'\leq q$ be the maximal positive integers such that
\[
\lambda_1=\dots=\lambda_{p'}\qquad\text{and}\qquad  \lambda_{n-q'+1}=\dots=\lambda_n. 
\]
Then 
\[
(\la-\beta)^+=\la-(\eps_{p'}-\eps_{n-q'+1}).
\]
The root $\gamma=\eps_{p'}-\eps_{n-q'+1}$, or the integers $p'$ and $q'$, contain the same information about $\la$ as the root systems $Q=R$ of \cite{EHW}. In fact, it is easy to see that $Q$ is the root system of the subalgebra $\frsu(p',q')$ of $\frg$ built on the first $p'$ and last $q'$ coordinates.

It was proved in \cite{PPST1} that the inequality \eqref{basic di su} is equivalent to
\eq
\label{basic di su 2}
\la_1-\la_n\leq -n+p'+q'.
\eeq

As in the other cases, the $K$-module  
$S(\frp^-)$ is multiplicity free, and its $K$-types are the Schmid modules $F_{-s}$ with lowest weights
$-s$, where 
\eq
\label{gen schmid su}
s=(b_1,\dots,b_p\,|\,0,\dots,0,-b_p,\dots,-b_1),
\eeq
where $b_1\geq\dots\geq b_p\geq 0$ are integers. 
All Schmid modules are nonnegative integer combinations of the basic Schmid modules
\[
s_i=(\underbrace{1,\dots,1}_i,0,\dots,0\bbar 0,\dots,0,\underbrace{-1,\dots,-1}_i)
\] 
for $i=1,\dots,p$.

As in the other cases, 
we start by examining the condition on $\la$ which ensures
that the Dirac inequality
\eq
\label{din1 su}
\|(\la-s_i)^++\rho\|^2\geq\|\la+\rho\|^2
\eeq
holds, where $i=1,\dots,\min(p',q')$. We already know that for $i=1$ this is just the basic inequality \eqref{basic di su}, or equivalently \eqref{basic di su 2}.

It was proved in \cite{PPST1} that \eqref{din1 su} is equivalent to
\eq
\label{din3 su}
\la_1-\la_n\leq -n+p'+q'-i+1.
\eeq  

We will see in Theorem \ref{thm su} that if $\la$ satisfies \eqref{din3 su} as an equality for some $i\in[1,\min(p',q')]$, then $L(\la)$ is unitary. All other $\la$ are handled by the following theorem which follows from the results of \cite{PPST1}.

\begin{thm}
\label{unit nonunit su}
As above, let 
\[
\la=(\underbrace{\la_1,\dots,\la_1}_{p'},\la_{p'+1},\dots,\la_p\bbar\la_{p+1},\dots,\la_{n-q'},\underbrace{\la_n,\dots,\la_n}_{q'}).
\]
Then:
\begin{enumerate}
\item If 
\[
\la_1-\la_n>-n+p'+q',
\]
then $L(\la)$ is not unitary.
\item If for some integer $i\in[1,\min(p',q')-1]$ 
\eq
\label{gaps su}
-n+p'+q'-i<\lambda_1-\la_n< -n+p'+q'-i+1,
\eeq
then $L(\la)$ is not unitary.
\item If
\eq
\label{cont su}
\lambda_1-\la_n< -n+p'+q'-\min(p',q')+1=-n+\max(p',q')+1,
\eeq
then $L(\lambda)=N(\lambda)$ and it is unitary. These $\la$ are in the continuous part of its line in the sense of \cite{EHW} ($\lambda,\la'$ are on the same line if $\la-\la'$ has all coordinates equal to each other).
\end{enumerate}
\end{thm}

\pf
(1) is immediate from our basic inequality \eqref{basic di su 2}.

(2) follows from 
Theorem 3.4.(1) 
of \cite{PPST1}, from Corollary \ref{cor unit nonunit}.(1), and from the obvious fact that a Schmid module $s$ of lower level than $s_{i +1}$ can have at most $i$ nonzero components on each side of the bar. 

(3) follows from 
Theorem 3.4.(2) 
of \cite{PPST1} and from  Corollary \ref{cor unit nonunit}.(2). 
\epf

To complete the classification, we prove the following theorem.

\begin{thm}\label{thm su}
Let $p',q',i$ be integers such that $1\leq p'\leq p$, $1\leq q'\leq q$ and $1\leq i\leq \min(p',q')$. Let $a_1\geq\ldots\geq a_{p-p'}\geq 1$ and $b_1\geq\ldots\geq b_{q-q'}\geq 1$ 
be integers. Let 
\eq\label{la1 su}
\la_1=-n+p'+q'-i+1
\eeq 
and let $L_\la$ be the irreducible module with highest weight 
\eq\label{unitary su}
\la=(\underbrace{\la_1,\ldots,\la_1}_{p'},
\la_1-a_{p-p'},\ldots,\la_1-a_1\bbar b_1,\dots,b_{q-q'},\underbrace{0,\dots,0}_{q'}).
\eeq
 Then $L(\la)$ is unitary.
\end{thm}

The cones of \cite{DES} are easily described in terms of \eqref{unitary su}. The vertex of the cone containing $\lambda$ is obtained by setting all $a_r$ and $b_s$  equal to 1, and all elements of the cone are obtained by taking all possible choices for the $a_r$ and $b_s$.

To prove Theorem \ref{thm su}, we first define the basic representations $L_k$, $k\in\bbZ$. The highest weight of $L_k$ is $\la_k$ given by
\eq\label{basic su}
\la_k=\left\{ \begin{matrix} (-1,\ldots,-1,-1-k\bbar 0,\dots,0), & k\geq 0 \\ (-1,\ldots,-1\bbar -k,\dots,0), & k\leq -1.\end{matrix}\right.
\eeq
(These modules are called ladder representations in \cite{DES}.)

\begin{lem}  
The modules $L_k$, $k\in\bbZ,$ are unitary.
\end{lem}

\pf 
We first describe an embedding of $\fru(p,q)_\bbC$ into $\frsp(2n,\bbR)_\bbC$ and then realize $L_k$ inside (the restriction of) the Weil representation. The embedding is well known; it is the special case $m=1$ of the dual pair $\fru(p,q)\times\fru(m)\subset\frsp(2m(p+q),\bbR)$. See for example \cite{NOT}, Section 3. We describe the embedding explicitly in terms of roots and root vectors as follows. 
This embedding is associated with the following Vogan diagram for $\frsp(2n,\bbR)$:

\[
\dynkin \whiteroot{} \link\whiteroot{}\dots\whiteroot{}\link\root{p}\link\whiteroot{}\dots \whiteroot{}
\llink< \root{}
\enddynkin
\]

To obtain a Vogan diagram for $\fru(p,q)$ we cross the last point in the above diagram. 
More explicitly, the simple roots of $\fru(p,q)_\bbC$ are 
\[
\eps_1-\eps_2,\dots,\eps_{p-1}-\eps_p,\eps_p+\eps_{p+1},\eps_{p+2}-\eps_{p+1},\dots,\eps_n-\eps_{n-1}.
\]
All positive roots of $\fru(p,q)_\bbC$ are thus
\begin{eqnarray*}
&\eps_i-\eps_j,&\quad 1\leq i<j\leq p;\\
&\eps_j-\eps_i,&\quad p+1\leq i<j\leq n; \\
&\eps_i+\eps_j,&\quad 1\leq i\leq p,\ p+1\leq j\leq n.
\end{eqnarray*}
The corresponding positive root vectors are
\begin{eqnarray*}
&\del_i z_j,&\quad 1\leq i<j\leq p;\\
&\del_j z_i,&\quad p+1\leq i<j\leq n; \\
&\del_i\del_j,&\quad 1\leq i\leq p,\ p+1\leq j\leq n.
\end{eqnarray*}
The Cartan subalgebra for $\fru(p,q)_\bbC$ is the same as the Cartan subalgebra for $\frsp(2n,\bbR)_\bbC$, but the basis is
\[
-z_1\del_1-\half,\dots,-z_p\del_p-\half,z_{p+1}\del_{p+1}+\half,\dots,z_n\del_n+\half.
\]
(The second group of basis elements has a change of sign because the second group of positive roots for $\fru(p,q)_\bbC$ consists of roots that are negative in the standard choice for $\frsp(2n,\bbR)_\bbC$.)

It is now clear that each of the vectors $z_p^k$, $k\in\bbZ_+$, is a highest weight vector for $\frsu(p,q)_\bbC$, of weight
\[
(-\half,\dots,-\half,-\half-k\bbar \half,\dots,\half),
\]
which becomes
\[
(-1,\dots,-1,-1-k\bbar 0,\dots,0)
\]
by subtracting $\half$ from all coordinates.
We now use Proposition \ref{dot} to generate an irreducible unitary highest weight module for $\frsu(p,q)$ by $z_p^k$. In this way we have constructed the basic representations $L_k$, $k\geq 0$.

To construct $L_k$, $k<0$, we set $j=-k>0$ and note that $z_{p+1}^j$ is a highest weight vector for $\frsu(p,q)_\bbC$, of weight
\[
(-\half,\dots,-\half\bbar\half+j,\half,\dots,\half)=(-1,\dots,-1\bbar j,0,\dots,0).
\]
By Proposition \ref{dot}, this vector generates an irreducible unitary highest weight module for $\frsu(p,q)_\bbC$. This module is clearly isomorphic to $L_k$.
\epf

{\it Proof of Theorem \ref{thm su}.} 
Let
\[
L_{a_1,\ldots, a_{p-p'};b_1,\dots,b_{q-q'}}=
 L_{a_{p-p'}}\bul(\ldots \bul(L_{a_1}\bul (L_{-b_{q-q'}}\bul(\ldots\bul L_{-b_1})\ldots))\ldots).
\]
and let
\[
(L_0)^{i-1}=\underbrace{L_0\bul(\dots \bul(L_0\bul L_0)\dots)}_{i-1}.
\]
We are going to show that
\eq\label{decomp su}
  L(\la)\cong L_{a_1,\ldots a_{p-p'};b_1,\dots,b_{q-q'}}\bul (L_0)^{i-1}.
\eeq
This will imply the theorem since PRV products of unitary modules are unitary. 
	
\begin{lem}\label{bul su}
Fix integers $r$ and $s$ such that $1\leq r\leq p'$ and $1\leq s\leq q'$. 
Let $ a_1\geq\ldots\geq a_r\geq 1$ and $b_1\geq\dots\geq b_s\geq 1$ be sequences of integers and consider the irreducible unitary highest weight module  
\[
L_{a_1,\ldots a_{r};b_1,\dots,b_{s}}=
 L_{a_{r}}\bul(\ldots \bul(L_{a_1}\bul (L_{-b_s}\bul(\ldots\bul L_{-b_1})\ldots))\ldots).
\]
Then $L_{a_1,\ldots, a_{r};b_1,\dots,b_{s}}$ has highest weight 
\[
\la_{a_1,\ldots,a_r;b_1\dots,b_s}=(\underbrace{-r-s,\ldots,-r-s}_{p-r},
-r-s-a_r,\ldots,-r-s-a_1\bbar b_1,\dots,b_s,\underbrace{0,\dots,0}_{q-s}).
\]
\end{lem}
 \pf
 We first use induction on $s$ to show that $L_{-b_s}\bul(\ldots\bul L_{-b_1})\ldots)$ has highest weight 
 \[
 (-s,\dots,-s\bbar b_1,\dots,b_s,0,\dots,0).
 \]
Then we use induction on $r$ to prove the statment of the lemma. The computations are very similar to the the proof of Lemma \ref{conus}, so we omit them.
\epf
 
We now get back to the proof of Theorem \ref{thm su}. 
Using Lemma \ref{bul su} for $r=p-p'$ and $s=q-q'$ we see that the highest weight of $ L_{a_1,\ldots, a_{p-p'};b_1,\dots,b_{q-q'}}$
is 
\[
(\underbrace{\alpha,\ldots,\alpha}_{p'},
\alpha-a_{p-p'},\ldots,\alpha-a_1\bbar b_1,\dots,b_{q-q'},\underbrace{0,\dots,0}_{q'}),
\]
where $\alpha=-n+p'+q'$. 
Since the lowest weight of the top $K$-type of $(L_0)^{i-1}$
is 
\[
(1-i,\dots,1-i\bbar 0,\dots,0),
\]
we see that $L_{a_1,\ldots a_{p-p'};b_1,\dots,b_{q-q'}}\bul (L_0)^{i-1}$ has highest weight $\la$. This proves \eqref{decomp su} and hence finishes the proof of Theorem \ref{thm su}.
\qed

\section{Classification for $\frso^*(2n)$}
\label{subsec di so*}

In this case $\frg=\frso(2n,\bbC)$ and $\frk=\frgl(n,\bbC)$. We use standard coordinates and standard positive root systems: the compact Cartan subalgebra is identified with $\bbC^n$, the positive compact roots are $\eps_i-\eps_j$, $i<j$, and the positive noncompact roots are $\eps_i+\eps_j$, $i<j$. So
\[
\rho=(n-1,n-2,\dots,1,0)
\]
and the highest (noncompact) root is 
\[
\beta=\eps_1+\eps_2.
\]
The highest weight $(\frg,K)$-modules have highest weights of the form
\[
\la=(\la_1,\la_2,\dots,\la_n),\qquad \la_1\geq\la_2\geq\dots\geq\la_n,\qquad \la_i-\la_j\in\bbZ,\ 1\leq i,j\leq n.
\]
The basic necessary condition for unitarity is, as before, the Dirac inequality
\eq
\label{basic di so}
\|(\la-\beta)^++\rho\|^2\geq\|\la+\rho\|^2.
\eeq
To make this inequality more precise, we as before write
\[
(\la-\beta)^+=\la-\gamma
\]
for a positive noncompact root $\gamma$. Then \eqref{basic di so} becomes 
\[
\|\la-\gamma+\rho\|^2\geq\|\la+\rho\|^2,
\]
and since $\|\gamma\|^2=2$, this is equivalent to
\[
\langle \la+\rho,\gamma\rangle\leq 1.
\]
There are two basic cases:
\smallskip

\noindent{\bf Case $\mathbf{1}$: $\mathbf{\la_1>\la_2}$.} Let $q\in[2,n]$ be such that $\la_2=\dots=\la_q$ and, in case $q<n$, $\la_q>\la_{q+1}$. Then $\gamma=\eps_1+\eps_q$, and since $\la_q=\la_2$, the basic inequality becomes
\eq 
\label{basic di so c1}
\la_1+\la_2\leq -2n+q+2.
\eeq

\noindent{\bf Case $\bold 2$: $\bold{\la_1=\la_2}$.} Let $p\in[2,n]$ be such that $\la_1=\la_2=\dots=\la_p$ and, in case $p<n$, $\la_p>\la_{p+1}$. Then $\gamma=\eps_{p-1}+\eps_p$, and since $\la_{p-1}=\la_p=\la_1$, the basic inequality becomes
\eq 
\label{basic di so c2}
\la_1\leq -n+p.
\eeq
\smallskip

Since the strongly orthogonal Harish-Chandra roots are
\[
\gamma_i=\eps_{2i-1}+\eps_{2i},\qquad i=1,2,\dots,[n/2],
\]
the basic Schmid $\frk$-submodules of $S(\frp^-)$ have lowest weights $-s_i$, where
\eq
\label{bas sch so}
s_i=(\underbrace{1,\dots,1}_{2i},0,\dots,0),\qquad i=1,2,\dots,[n/2].
\eeq
(Of course, $s_1=\gamma_1=\beta$.) Moreover, all irreducible $\frk$-submodules of $S(\frp^-)$ have lowest weights $-s$, where
\eq
\label{gen sch so}
s=(b_1,b_1,b_2,b_2,\dots,b_j,b_j,0,\dots,0)
\eeq
for some $j$, $1\leq j\leq [n/2]$, and some positive integers $b_1\geq b_2\geq\dots\geq b_j$. 

\smallskip

We first describe the unitary modules in Case 1. In this case, there is only one discrete point on each line, so this point is the limit of the continuous part of the line. The following result follows from 
Theorem 3.2 
of \cite{PPST1} and from Corollary \ref{cor unit nonunit} by the same arguments we used to prove Theorem \ref{unit nonunit sp} and Theorem \ref{unit nonunit su}.

\begin{thm}
\label{thm so c1} 
Suppose $\la$ is as in Case 1, i.e., there is $q\in [2,n]$ such that
\[
\la_1>\la_2=\dots=\la_q,
\]
and $\la_q>\la_{q+1}$ if $q<n$. Then $L(\la)$ is unitary if and only if $\la$ satisfies the basic Dirac inequality \eqref{basic di so c1}. 

If $\la$ satisfies \eqref{basic di so c1} strictly, then $L(\la)$ is equal to the full generalized Verma module $N(\lambda)$, i.e., $N(\lambda)=L(\la)$ is irreducible and unitary. Such a module is in the continuous part of its line in the sense of \cite{EHW}.
\end{thm}

We now turn to Case 2, i.e., 
\[
\la=(\underbrace{\la_1,\dots,\la_1}_p,\la_{p+1},\dots,\la_n)
\]
for some $p\in[2,n]$, with $\la_1>\la_{p+1}$ if $p<n$. Besides the basic inequality \eqref{basic di so c2}, we also examine when
\eq
\label{di si}
\|(\la-s_i)^++\rho\|^2\geq \|\la+\rho\|^2
\eeq
for a basic Schmid module $s_i$ with $2i\leq p$. It was proved in \cite{PPST1} that \eqref{di si} is equivalent to
\eq
\label{disc pts so}
\la_1\leq -n+p-i+1.
\eeq
In particular, for $i=1$ this is the basic inequality \eqref{basic di so c2}.

We will show in Theorem \ref{thm so*} that if $\la$ belonging to Case 2 is such that \eqref{disc pts so} is an equality for some integer $i\in [1,[\frac{p}{2}]]$ then $L(\la)$ is unitary. The other $\la$ belonging to Case 2 are handled by the following theorem which follows from 
Theorem 3.3 
of \cite{PPST1} and from Corollary \ref{cor unit nonunit} by the same arguments we used to prove Theorem \ref{unit nonunit sp} and Theorem \ref{unit nonunit su}.

\begin{thm}
\label{unit nonunit so}
Let $\la$ be in Case 2, i.e.,
\[
\la=(\underbrace{\la_1,\dots,\la_1}_p,\la_{p+1},\dots,\la_n)
\]
for some $p\in[2,n]$, with $\la_1>\la_{p+1}$ if $p<n$. Then:
\begin{enumerate}
\item If 
\[
\la_1>-n+p,
\]
then $L(\la)$ is not unitary.
\item If for some integer $i\in[1,[\frac{p}{2}]-1]$ 
\eq
\label{gaps so}
-n+p-i<\lambda_1< -n+p-i+1,
\eeq
then $L(\la)$ is not unitary.
\item If
\eq
\label{cont so}
\lambda_1< -n+p-\left[\frac{p}{2}\right]+1=-n+\left[\frac{p+1}{2}\right]+1,
\eeq
then $L(\lambda)=N(\lambda)$ and it is unitary. (Such a $\la$ is in the continuous part of its line in the sense of \cite{EHW}.) 
\end{enumerate}
\end{thm}

Note that the classification is already obtained for Case 1 by Theorem \ref{thm so c1}, since every point that can be unitary is either in the continuous part of its line and hence unitary by Theorem \ref{thm so c1}, or at the end of the continuous part of its line and therefore unitary by Lemma \ref{continuity} (see also \cite[Proposition 3.1.c]{EHW}).

We can thus concentrate on Case 2, and to
complete the classification, we prove the following theorem.

\begin{thm}\label{thm so*}
Let $p$ and $i$ be integers such that $2\leq p\leq n$ and $1\leq i\leq \left[\frac{p}{2}\right]$. Let $a_1\geq\dots\geq a_{n-p}\geq 1$ be integers. 
Let 
\eq\label{la1 so*}
\la_1=-n+p-i+1
\eeq 
and let $L_\lambda$ be the irreducible module with highest weight 
\eq\label{unitary so*}
\la=(\underbrace{\la_1,\dots,\la_1}_p,\la_1-a_{n-p},\dots,\la_1-a_1).
\eeq
 Then $L_\la$ is unitary.
\end{thm}

The cones of \cite{DES} are easily described in terms of \eqref{unitary so*}. The vertex of the cone containing $\lambda$ is obtained by setting all $a_k$ equal to 1, and all elements of the cone are obtained by taking all possible choices for the $a_k$.

For the proof of Theorem \ref{thm so*} we first describe the basic representations:

\eq\label{basic so*}
L_{k}:=L_{\la_k};\;\la_k=(-1,\ldots,-1,-1-k);\;k\geq 0.
\eeq

\begin{lem}  
The modules $L_{k}$, $k\geq 0,$ are unitary.
\end{lem}

\pf We first describe a (well known) embedding of $\frso^*(2n)_\bbC=\frso(2n,\bbC)$ into $\frsu(n,n)_\bbC=\frsl(2n,\bbC)$ and then construct the basic representations for $\frso^*(2n)$ inside the restrictions of the basic representations for $\frsu(n,n)$.

We consider $\frso(2n,\bbC)$ corresponding to the symmetric bilinear form given by the antidiagonal matrix
\[
S=\begin{pmatrix} 0&0&\dots&0&1\\
0&0&\dots&1&0\\
\dots&\dots&\dots&\dots&\dots\\
0&1&\dots&0&0\\
1&0&\dots&0&0
\end{pmatrix}
\]
Then $\frso(2n,\bbC)$ consists of the matrices that are skew symmetric with respect to the antidiagonal. These are the $2n\times 2n$ matrices $X$ that satisfy the condition 
\[
X^{t'}=-X,
\]
where $X^{t'}_{ij}=X_{2n+1-j\,2n+1-i}$. The standard Cartan decomposition of $\frso(2n,\bbC)$ corresponding to the real form $\frso^*(2n)$ is given by
\begin{eqnarray*}
&\frk&=\left\{\begin{pmatrix} A & 0\\ 0 &-A^{t'}\end{pmatrix}\bbar A\in\frgl(n,\bbC)\right\};\\
&\frp^+&=\left\{\begin{pmatrix} 0 & B\\ 0 &0\end{pmatrix}\bbar B\in\frgl(n,\bbC),B^{t'}=-B\right\};\\
&\frp^-&=\left\{\begin{pmatrix} 0 & 0\\ C &0\end{pmatrix}\bbar C\in\frgl(n,\bbC),C^{t'}=-C\right\}.
\end{eqnarray*}
On the other hand, the standard Cartan decomposition of $\frsl(2n,\bbC)$ corresponding to the real form $\frsu(n,n)$ has $\frk$ consisting of block diagonal matrices, $\frp^+$ consisting of block upper triangular matrices and $\frp^-$ consisting of block lower triangular matrices. So we see that the embedding of $\frso(2n,\bbC)$ into $\frsl(2n,\bbC)$ is compatible with standard Cartan decompositions.

Let $\frh$ be the Cartan subalgebra of $\frsl(2n,\bbC)$ consisting of traceless diagonal matrices. The intersection \[
\frt=\frh\cap\frso(2n,\bbC)
\] 
is a Cartan subalgebra of $\frso(2n,\bbC)$; it consists of diagonal matrices with diagonal entries
\[
a_1,a_2,\dots,a_{n-1},a_n,-a_n,-a_{n-1},\dots,-a_2,-a_1.
\]
The standard choice for positive roots of $\frsl(2n,\bbC)$ with respect to $\frh$ is $\eps_i-\eps_j$, $1\leq i<j\leq 2n$, while the standard choice for positive roots of $\frso(2n,\bbC)$ with respect to $\frt$ is $\eps_i'\pm\eps_j'$, $1\leq i<j\leq n$. Here $\eps_k'$ denotes the restriction of $\eps_k$ to $\frt$.

For any $1\leq i<j\leq n$, $\eps_i-\eps_j$ and $\eps_{2n+1-j}-\eps_{2n+1-i}$ both have restrictions to $\frt$ equal to $\eps_i'-\eps_j'$. The root vector corresponding to the positive root $\eps_i'-\eps_j'$ of $\frso(2n,\bbC)$ is
$E_{ij}-E_{2n+1-j\,2n+1-i}$. 

For any $1\leq i<j\leq n$, $\eps_i-\eps_{2n+1-j}$ and $\eps_j-\eps_{2n+1-i}$ both have restrictions to $\frt$ equal to $\eps_i'+\eps_j'$. The root vector corresponding to the positive root $\eps_i'+\eps_j'$ of $\frso(2n,\bbC)$ is
$E_{i\,2n+1-j}-E_{j\,2n+1-i}$. 

So we see that the embedding of $\frso(2n,\bbC)$ into $\frsl(2n,\bbC)$ is compatible with standard choices of positive roots and that positive root vectors for $\frso(2n,\bbC)$ are linear combinations of positive root vectors for $\frsl(2n,\bbC)$.

Let now $\tilde L_k$ be the basic representation of $\frsu(n,n)_\bbC=\frsl(2n,\bbC)$ with highest weight
\[
\tilde\la_k=(-1,\dots,-1,-1-k\bbar 0,\dots,0),\qquad k\geq 0,
\]
and let $w$ denote its highest weight vector. Because of the above compatibility, $w$ is a highest weight vector also for $\frso(2n,\bbC)$. Its weight is the restriction of $\tilde\la_k$, which is equal to
\[
\la_k=(-1,\dots,-1,-1-k).
\]
By Proposition \ref{dot}, $w$ generates an irreducible unitary highest weight module for $\frso^*(2n)$ with highest weight $\la_k$. Thus we have found a unitary realization of $L_k$.
\epf

{\it Proof of Theorem \ref{thm so*}.} 
Let $L_0^{i-1}=
 \underbrace{L_{0}\bul (L_0\bul(\ldots \bul(L_{0}\bul L_0)\ldots))}_{i-1}$ and let
 \[
 L_{a_1,\dots,a_{n-p}}=\left( L_{a_{n-p}}\bul (L_{a_{n-p-1}}\bul(\ldots \bul(L_{a_2}\bul L_{a_1})\ldots))\right).
\]
We are going to show that
\eq\label{decomp so*}
L(\la)\cong  L_{a_1,\dots,a_{n-p}}\bul (L_0)^{i-1}.
\eeq
This will imply the theorem since PRV products of unitary modules are unitary. 
	
\begin{lem}\label{bul so*}
Fix an integer $j$ such that $1\leq j\leq n$. 
Let $ a_1\geq\ldots\geq a_j\geq 1$ be integers.
Then the module 
\[
L_{a_1,\dots,a_j}= \left( L_{a_j}\bul (L_{a_{j-1}}\bul(\ldots \bul(L_{a_2}\bul L_{a_1})\ldots))\right)
\]
has highest weight
\[
\la_{a_1,\dots,a_j}=(\underbrace{-j,\ldots,-j}_{n-j},-j-a_j,\ldots,-j-a_1).
\]
\end{lem}
 
 \pf
 We use induction on $j$. The case $j=1$  is trivial.
 
 Suppose that the statement holds for the sequence $a_1\geq\ldots\geq a_{j-1}$, i.e., that $L_{a_1,\ldots a_{j-1}}$ has highest weight
 \[
 (\underbrace{-j+1,\ldots,-j+1}_{n-j+1},
-j+1-a_{j-1},\ldots,-j+1-a_1).
\]
Consider the  product $L_{a_j}\bul L_{a_1,\ldots, a_{j-1}}.$ The lowest weight of the top K-type
 of $L_{a_1,\ldots, a_{j-1}}$ is equal to $(-j+1-a_1,\ldots,-j+1-a_{j-1},-j+1,\ldots,-j+1)$, so $L_{a_j}\bul L_{a_1,\ldots, a_{j-1}}$  
 has highest weight 
 \[
 (-j-a_1,\ldots,-j-a_{j-1},
 -j,\ldots,-j,-j-a_j)^+
\] 
 which is equal to $\la_{a_1,\ldots, a_j}.$ 
 \epf
 
We now get back to the proof of Theorem \ref{thm so*}. From   Lemma \ref{bul so*} applied for $j=n-p$, we see that   
the highest weight of $L_{a_1,\ldots, a_{n-p}}$ is 
$$
(\underbrace{-n+p,\ldots,-n+p}_{p},
-n+p-a_{n-p},\ldots,-n+p-a_1).
$$
Since the lowest weight of the top $K$-type of $(L_0)^{i-1}$
is 
$$
(-i+1,\ldots,-i+1),
$$
we see that $L_{a_1,\ldots a_{n-p}}\bul (L_0)^{i-1}$ has highest weight $\la$. This proves \eqref{decomp so*} and hence finishes the proof of Theorem \ref{thm so*}.
\qed

\section{Unitary modules for $\frsp(2n,\bbR)$ with fixed integral or half integral infinitesimal character}
\label{sec:unitreg}

Recall that integrality of a weight $\mu$ means that in standard coordinates, the components of $\mu$ are all integers. We will also consider the cases when all coordinates are half integers (elements of $\half+\bbZ$). An infinitesimal character $\chi_\La$ is integral (respectively half integral) if the weight $\La$ corresponding to it by the Harish-Chandra homomorphism is integral (respectively half integral). In the following we identify $\chi_\La$ with the weight $\La$ (modulo the Weyl group action). The reason for considering these cases is the fact that all discrete unitary points in the classification are either integral or half integral, and these are the most interesting cases. 

By definition, the parameter of the irreducible highest weight $(\frg,K)$-module $L(\la)$ is $\La=\la+\rho$. It is a particular representative of the infinitesimal character of $L(\la)$. If $\la$ is $\frg$-dominant, then any irreducible highest weight $(\frg,K)$-module with the same infinitesimal character as $L(\la)$ is of the form $L(w(\la+\rho)-\rho)$ for some $w$ in
\[
W^1=\{w\in W_\frg\,\big|\, w\rho \text{ is $\frk$-dominant}\}.
\]
The corresponding parameter is $w(\la+\rho)$. These parameters, or the corresponding Weyl group elements $w$, can be organized into a Hasse diagram with respect to the Bruhat order. We will show examples of that later in this section.

We will denote by $\La^{\dom}$ the $\frg$-dominant representative of the infinitesimal character $\La$. All parameters $\La$ corresponding to highest weight $(\frg,K)$-modules must be dominant regular for $\frk$. This means that $\La^{\dom}$ can have each coordinate repeated at most twice, and can contain zero at most once. If a certain coordinate $z$ is repeated in $\La^{\dom}$, then all parameters $\La$ must contain $z$ and $-z$.

Recall that in the integral or half-integral case $L(\la)$ is unitary if and only if
\[
\la_1\leq -n+\frac{q+r}{2},
\]
where $q\leq r$ are integers such that
\[
\la=(\underbrace{\la_1,\dots,\la_1}_q,\underbrace{\la_1-1,\dots,\la_1-1}_{r-q},\la_{r+1},\dots,\la_n).
\]
In terms of the corresponding parameter
\[
\La=\la+\rho=(\la_1+n,\la_2+n-1,\dots,\la_n+1),
\]
this means the following:
\begin{enumerate}
\item $\La_1,\dots,\La_q$ is a string (a sequence of integers descending by 1);
\item if $r>q$, then $\La_{q+1},\dots,\La_r$ is another string with the gap between $\La_q$ and $\La_{q+1}$ equal to 2; 
\item the gap between $\La_r$ and $\La_{r+1}$ is at least 2 if $r>q$ and at least 3 if $r=q$;
\item the unitarity condition holds:
\eq
\label{unit La}
\La_1\leq\frac{q+r}{2}.
\eeq
\end{enumerate}
If the module $L(\la)$ is unitary, we will say that the corresponding $\La$ is a unitary parameter.

\begin{rem} 
\label{neg}
{\rm
It is clear from the above that if $\La_1\leq 0$, then $\La$ is unitary. There is at most one such parameter for a fixed $\La^{\dom}$. This means that, whenever necessary, we can assume that $\La_1>0$.

Note that $\La_1\leq 0$ can only happen when there are no repeated coordinates in $\La^{\dom}$, so $\La$ is either regular, or has $\La_1=0$ as the only singularity.
}
\end{rem}

We now want to describe a criterion for a parameter $\La$ to be unitary that will enable us to write down the list of all unitary $\La$ for a fixed $\La^{\dom}$.

Let $x\geq 0$ be the maximal number such that the coordinates $(\La^{\dom})_i$ of $\La^{\dom}$ include
\eq
\label{x int}
0,1,1,\dots,x,x\qquad \text{if all } (\La^{\dom})_i\in\bbZ,\quad\text{or}
\eeq
\eq
\label{x hint}
\half,\half,\dots,x,x\qquad \text{if all } (\La^{\dom})_i\in\half+\bbZ.
\eeq
Then each parameter $\La\in W^1\La^{\dom}$ of a highest weight $(\frg,K)$-module contains the string
\eq
\label{x string}
x,x-1,\dots,-x.
\eeq
It is also possible that $x$ as above does not exist, i.e., that $\La^{\dom}$ does not contain 0 (if the coordinates are integers) or two copies of $\half$ (if the coordinates are half integers). Note that if the string $x,\dots,-x$ is a substring of a string $c_1,\dots,c_k$, then, by maximality of $x$, the string $x,\dots,-x$ is either the leftmost part or the rightmost part of the string $c_1,\dots,c_k$. In other words, either $c_1=x$ or $c_k=-x$.

\begin{thm}
\label{unit crit}
(1) 
Let $\La^{\dom}$ be such that $x$ as in \eqref{x int} or \eqref{x hint} exists, and suppose $x$ is the maximal such number. Let $\La$ be a parameter; it has to contain the string $x,x-1,\dots,-x$ (if $x=0$, then this string consists of just 0).
\begin{enumerate}
\item[(1a)] If the string $x,\dots,-x$ occurs inside the string $\La_1,\dots,\La_q$, then $\La$ is unitary.
\item[(1b)] If $r>q$ and the string $x,\dots,-x$ occurs inside the string $\La_{q+1},\dots,\La_r$, then $\La$ is unitary if and only if $x,\dots,-x$ is the leftmost part of, and not all of, $\La_{q+1},\dots,\La_r$. Equivalently, $\La_r\leq -x-1$, and then necessarily $\La_{q+1}=x$. (In particular, if $\La$ is unitary, then $r\geq q+2$.)
\item[(1c)] If the string $x,\dots,-x$ occurs inside the group of coordinates $\La_{r+1},\dots,\La_n$, then $\La$ is not unitary.
\end{enumerate}

(2) Let $\La^{\dom}$ be such that $x$ as in \eqref{x int} or \eqref{x hint} does not exist, i.e., $\La^{\dom}$ does not have a coordinate equal to 0, or two coordinates equal to $\half$.
\begin{enumerate}
\item[(2a)] If $\La$ is such that $r=q$, then $\La$ is unitary if and only if $\La_q\leq 1$.
\item[(2b)] If $\La$ is such that $r>q$, then $\La$ is unitary if and only if $\La_q\leq \frac{3}{2}$.
\end{enumerate}
\end{thm}
\pf
(1a) Leq $a,b\geq 0$ be integers such that
\[
\La_1,\dots,\La_q=x+a,\dots,x,\dots,-x,\dots,-x-b.
\]
(By maximality of $x$, at least one of $a,b$ must be 0.) Since the length of the string $x,\dots,-x$ is $2x+1$, we see that
$q= 2x+1+a+b$. Therefore
\[
\La_1=x+a< q\leq \frac{q+r}{2},
\]
so \eqref{unit La} holds and $\La$ is unitary.
\smallskip

(1b) 
Let $a,b\geq 0$ be integers such that
$$
\Lambda_{q+1},\ldots,\Lambda_r=x+a,\ldots,-x-b.
$$
Then by maximality of the string $x,\dots,-x$, either $a$ or $b$ must be 0.

If $a=0,$ then $r-q=2x+1+b, $ hence
$$
\frac{q+r}{2}=x+q+\frac{b+1}{2}.
$$
On the other hand, $\Lambda_{q+1}=x$ implies $\Lambda_1=x+q+1.$
Hence $\Lambda$ is unitary if and only if $b\geq 1,$ which is equivalent to $\La_r\leq -x-1$.

If $b=0$, then $r-q=2x+1+a, $ hence
$$
\frac{q+r}{2}=x+q+\frac{a+1}{2}.
$$
On the other hand, $\Lambda_{q+1}=x+a$ implies $\Lambda_1=x+a+q+1.$ Since $a\geq 0$, it follows that $\La_1>\frac{q+r}{2}$, so $\Lambda$ is not unitary.
\smallskip

(1c) The assumption implies that $\La_r\geq x+2$, and this implies that $\La_1\geq x+r+1$. In particular, 
\[
\La_1>r\geq \frac{q+r}{2},
\]
so \eqref{unit La} fails and $\La$ is not unitary.
\smallskip

(2a) Since $r=q$, the unitarity condition \eqref{unit La} is $\La_1\leq q$, which is equivalent to $\La_q\leq 1$.
\smallskip

(2b) If $\La_q\leq\frac{3}{2}$ then 
\[
\La_1\leq q+\half=\frac{2q+1}{2}\leq \frac{q+r}{2},
\] 
so \eqref{unit La} holds and $\La$ is unitary.

Conversely, if $\La_q\geq 2$, then $\La_{q+1}\geq 0$, so $\La_{q+1}>0$ since 0 is not a coordinate of $\La$. It follows that $\La_r>0$ since $\La$ does not contain 0 or the string $\half,-\half$. Therefore
\[
\La_1=\La_r+r>r\geq \frac{q+r}{2},
\]
so \eqref{unit La} fails and $\La$ is not unitary.
\epf

\begin{ex}
\label{ex sp 1}
{\rm
Let
\[
\La^{\dom}=(7,5,4,4,3,2,2,1,1,0).
\]
We want to use Theorem \ref{unit crit} to write the list of all unitary parameters conjugate to $\La^{\dom}$.

The number $x$ from Theorem \ref{unit crit} equals 2, and any parameter $\La$ of a highest weight $(\frg,K)$ module must contain the string $2,1,0,-1,-2$. Suppose first that this string is inside the initial string $\La_1,\dots,\La_q$ of $\La$, i.e., we are in the case (1a) of Theorem \ref{unit crit}. Then, since 4 is a repeated coordinate, $\La_1,\dots,\La_q$ must contain the string
\[
4,3,2,1,0,-1,-2.
\]
If $\La$ would have 7 for a coordinate, then the initial string would consists of 7 only, which is a contradiction. So $\La$ must have -7 for a coordinate. It remains to decide whether $\La$ has coordinate 5 or -5; both choices lead to unitary modules. So we get two unitary parameters 
\begin{gather*}
\La=(4,3,2,1,0,-1,-2,-4,-5,-7);\\
\La=(5,4,3,2,1,0,-1,-2,-4,-7).
\end{gather*}
Suppose now that the string $2,1,0,-1,-2$ is inside the second string $\La_{q+1},\dots,\La_r$ of $\La$, i.e., we are in the case (1b) of Theorem \ref{unit crit}. By (1b) of Theorem \ref{unit crit}, the string $\La_{q+1},\dots,\La_r$ must start with $x=2$, and contain the string $2,1,0,-1,-2,-3$. Moreover, this string continues with the repeated coordinate -4, while $\La_q$ must be equal to 4. As before, we must have -7 at the end. We again have a choice whether $\La$ has 5 or -5 as a coordinate, which leads us to two more unitary parameters
\begin{gather*}
\La=(4,2,1,0,-1,-2,-3,-4,-5,-7);\\
\La=(5,4,2,1,0,-1,-2,-3,-4,-7).
\end{gather*}
By Theorem \ref{unit crit}, these four parameters exhaust the list of unitary parameters $\La$ conjugate to $\La^{\dom}$.
}
\end{ex}
\begin{ex}
\label{ex sp 2}
{\rm
Let
\[
\La^{\dom}=(\frac{11}{2},\frac{9}{2},\frac{7}{2},\frac{7}{2},\frac{5}{2},\frac{3}{2},\frac{3}{2},\frac{1}{2}).
\]
We again want to use Theorem \ref{unit crit} to write the list of all unitary parameters conjugate to $\La^{\dom}$.

In this case, there is no $x$ as in Theorem \ref{unit crit}, i.e, the case (2) of the Theorem applies. We know that for $\La$ to be unitary, the string $\La_1,\dots,\La_q$ must end either with $\frac{3}{2}$ or with $\half$. Namely, we cannot have negative $\La_q$, because $\La_1\geq\frac{7}{2}$ since $\frac{7}{2}$ is a repeated coordinate, and thus also $\La_q$ must be positive since $\half$ is not a repeated coordinate and a string which starts positive and ends negative would have to contain $\half,-\half$.

If $\La_q=\frac{3}{2}$, then $\La_{q+1}=-\half$, so we see that $r>q$, and by Theorem \ref{unit crit} (2b), any such $\La$ will be unitary. It remains to see what choices we have for the string $\La_1,\dots,\La_q$. Since $\frac{7}{2}$ is a repeated coordinate, $\La_1,\dots,\La_q$ must include $\frac{7}{2},\frac{5}{2},\frac{3}{2}$, and we can freely choose to start the string with $\frac{7}{2}$, with $\frac{9}{2}$, or with $\frac{11}{2}$. This leads us to three unitary $\La$:
\begin{gather*}
\La=(\frac{7}{2},\frac{5}{2},\frac{3}{2},-\half,-\frac{3}{2},-\frac{7}{2},-\frac{9}{2},-\frac{11}{2});\\
\La=(\frac{9}{2},\frac{7}{2},\frac{5}{2},\frac{3}{2},-\half,-\frac{3}{2},-\frac{7}{2},-\frac{11}{2});\\
\La=(\frac{11}{2},\frac{9}{2},\frac{7}{2},\frac{5}{2},\frac{3}{2},-\half,-\frac{3}{2},-\frac{7}{2}).
\end{gather*}
If $\La_q=\frac{1}{2}$, analogous reasoning leads to three more unitary parameters:
\begin{gather*}
\La=(\frac{7}{2},\frac{5}{2},\frac{3}{2},\half,-\frac{3}{2},-\frac{7}{2},-\frac{9}{2},-\frac{11}{2});\\
\La=(\frac{9}{2},\frac{7}{2},\frac{5}{2},\frac{3}{2},\half,-\frac{3}{2},-\frac{7}{2},-\frac{11}{2});\\
\La=(\frac{11}{2},\frac{9}{2},\frac{7}{2},\frac{5}{2},\frac{3}{2},\half,-\frac{3}{2},-\frac{7}{2}).
\end{gather*}
By Theorem \ref{unit crit}, these six parameters exhaust the list of unitary parameters $\La$ conjugate to $\La^{\dom}$.
}
\end{ex}

We now generalize the arguments from Example \ref{ex sp 1} and Example \ref{ex sp 2} and use Theorem \ref{unit crit} to write down the list of unitary parameters $\La$ conjugate to a given $\La^{\dom}$.

\begin{cor}
\label{cor unit crit}
(1) Let $\La^{\dom}$ be such that (maximal) $x$ as in Theorem \ref{unit crit} exists.

(1a) Let $u\geq 0$ be the maximal integer such that $x+1,\dots,x+u$ are coordinates of $\La^{\dom}$ ($u=0$ when $x+1$ is not a coordinate of $\La^{\dom}$). Among these coordinates, let $x+v$, $0\leq v\leq u$, be the largest repeated coordinate (note that $v\neq 1$, since $x+1$ can not be repeated by maximality of $x$). Then the unitary $\La$ with $x,\dots,-x$ inside $\La_1,\dots,\La_q$ start with the string
\[
x+t,\dots,x,\dots,-x
\]
for some $t$, $v\leq t\leq u$, and have all other coordinates negative. (In particular, all coordinates of $\La^{\dom}$ greater than $x+u$ must be non-repeated.)

(1b) Assume that $x+1$ and $x+2$ are coordinates of $\La^{\dom}$, with $x+1$ non-repeated. Let $u\geq 2$ be the maximal integer such that $x+1,\dots,x+u$ are coordinates of $\La^{\dom}$, and let $v$, $0\leq v\leq u$, be such that $x+v$ is the maximal repeated coordinate among $x,\dots,x+u$ (note that $v\neq 1$). Then the unitary $\La$ with $x,\dots,-x$ inside $\La_{q+1},\dots,\La_r$ start with
\[
x+t,\dots,x+2,x,\dots,-x,-x-1
\]
for some $t$, $t\geq 2$, $v\leq t\leq u$, and have all other coordinates negative. (In particular, all coordinates of $\La^{\dom}$ greater than $x+u$ must be non-repeated.)

(2) Let $\La^{\dom}$ be such that $x$ as in Theorem \ref{unit crit} does not exist, i.e., $\La^{\dom}$ does not have a coordinate equal to 0 or two coordinates equal to $\half$. Then $\La$ is unitary if and only if one of the following conditions holds:
\begin{enumerate}
\item[(2a)] all coordinates of $\La$ are negative;
\item[(2b)] $\La$ starts with 
$q-\half,q-\frac{3}{2},\dots,\half$
for some $q\leq n$ and has other coordinates $\leq -\frac{3}{2}$; 
\item[(2c)] $\La$ starts with $q,q-1,\dots,1$
for some $q\leq n$ and has other coordinates $\leq -1$;
\item[(2d)]  $\La$ starts with $q+\half,q-\half,\dots,\frac{3}{2},-\half$
for some $q\leq n-1$ and has other coordinates $\leq -\frac{3}{2}$. 
\end{enumerate}
(Of course, in each case the coordinates involved must appear in $\La^{\dom}$.)
\end{cor}

\pf
(1a) We determine the $\La$ such that the string $x,\dots,-x$ is  in the first string of coordinates $\La_1,\dots,\La_q$; by Theorem \ref{unit crit} (1a), all such $\La$ are unitary. Note that some of the coordinates $x+2,\dots,x+u$ may be repeated, but $x+1$, if it occurs (i.e., if $u>0$) is not repeated, by maximality of $x$. 

Since $x+v$ is repeated, and since $x,\dots,-x$ is inside the initial string, $\La$ must contain the string
\eq
\label{string 1}
x+v,\dots,x+1,x,\dots,-x.
\eeq
If $v=u$, then $\La$ starts with \eqref{string 1} and continues with the negative coordinates.
If $v<u$, we can either start $\La$ with \eqref{string 1}, or continue to the left by $x+v+1,\dots,x+t$, where $v+1\leq t\leq u$ is chosen arbitrarily. 
\smallskip

(1b) We know from Theorem \ref{unit crit} (1b) that if $x,\dots,-x$ is inside $\La_{q+1},\dots,\La_r$, then $\La$ is unitary if and only if the string $\La_{q+1},\dots,\La_r$ starts with $x,\dots,-x,-x-1$; moreover $\La_q$ must be equal to $x+2$. In particular, both $x+1$ and $x+2$ must be coordinates of $\La^{\dom}$, and $x+1$ can not be a repeated coordinate.

Since $x+v$ is a repeated coordinate, and since $\La_1,\dots,\La_q=x+2$ is a string, $\La$ must contain the coordinates
\eq
\label{string 2}
x+v,\dots,x+2,x,\dots,-x,-x-1.
\eeq
If $v=u$, then $\La$ starts with \eqref{string 2} and continues with the negative coordinates. If $v<u$, we can either start $\La$ with \eqref{string 2} (so $t=v$), or continue to the left by $x+v+1,\dots,x+t$, where $t$ such that $t\geq 2$ and $v+1\leq t\leq u$, is chosen arbitrarily. 
\smallskip

(2) 
If all coordinates are negative, then $\La$ is unitary by
 Theorem \ref{unit crit}, or by Remark \ref{neg}, so we can assume $\La_1>0$. Then also $\La_q>0$, since $\La_1,\dots,\La_q$ is a string, and $\La$ does not have a coordinate equal to 0, or two coordinates equal to $\half,-\half$.
 
 By Theorem \ref{unit crit}, $\La$ is unitary if and only if:
 
either $\La_q=\half$, so $\La$ starts with
$q-\half,q-\frac{3}{2},\dots,\half$
for some $q\leq n$ and has other coordinates $\leq -\frac{3}{2}$; 
\smallskip

or $\La_q=1$, so $\La$ starts with
$q,q-1,\dots,1$
for some $q\leq n$ and has other coordinates $\leq -1$;
\smallskip

or $\La_q=\frac{3}{2}$ and $\La_{q+1}=-\half$ (so that $r>q$). So $\La$ starts with
$q+\half,q-\half,\dots,\frac{3}{2},-\half$
for some $q\leq n-1$ and has other coordinates $\leq -\frac{3}{2}$.
\epf

We now consider the special case
\[
\La^{\dom}=\rho=(n,n-1,\dots,1).
\]
All $\frk$-dominant $W$-conjugates of $\rho$ are of the form
\[
\La=(i_1,\dots,i_k,-j_m,\dots,-j_1)
\]
for various choices of integers $k,m\in[0,n] $ and $i_1>\dots>i_k$, $j_1>\dots>j_m$, such that $k+m=n$ and
\[
\{i_1,\dots,i_k\}\cup\{j_1,\dots,j_m\}=\{1,\dots,n\}.
\]
Following \cite{H2}, we express these parameters in terms of Young diagrams
with at most $n$ rows, with each row of length at most $n+1$, and with the additional condition that the diagrams can be built from hooks. The hooks are Young diagrams with row lengths
\[
k+1,\underbrace{1,\dots,1}_{k-1}\qquad\text{for some}\quad k\in[1,n].
\]
A Young diagram $Y$ is built from hooks if it is either a hook, or the first row and column of $Y$ form a hook, and after removing that hook the resulting Young diagram can still be built from hooks.

If the Young diagram $Y$ has row lengths $y_1\geq\dots \geq y_n$,  some of them possibly zero, we will write $Y=(y_1,\dots,y_n)$. The corresponding parameter $\La$ and highest weight $\la$ are
\[
\La=\rho-(y_n,\dots,y_1);\qquad \la=\La-\rho= -(y_n,\dots,y_1).
\]
\medskip

\begin{ex}
\label{ex yd sp}
{\rm Let $n=3$. The Hasse diagram of $\rho$ is

\medskip

{\footnotesize
\[
\begin{CD}
(3,2,1)@<<<(3,2,-1)@.@. \\
@. @AAA @.@. \\
@.  (3,1,-2)@<<<(3,-1,-2) @. \\
@.@AAA @AAA @.\\
@. (2,1,-3)@<<<(2,-1,-3) @. \\
@.@.@AAA@.\\
@.@.(1,-2,-3)@<<<(-1,-2,-3)
\end{CD}
\]
}
\medskip

\noindent with arrows pointing to larger elements in Bruhat order. The corresponding Young diagrams are

\medskip

{\tiny
\[
\begin{CD}
\emptyset@<<<\ydiagram{2}@.@. \\
@. @AAA @.@. \\
@.  \ydiagram{3,1}@<<<\ydiagram{3,3} @. \\
@.@AAA @AAA @.\\
@. \ydiagram{4,1,1}@<<<\ydiagram{4,3,1} @. \\
@.@.@AAA@.\\
@.@.\ydiagram{4,4,2}@<<<\ydiagram{4,4,4}
\end{CD}
\]
}
\medskip

}
\end{ex}

\medskip

Using Corollary \ref{cor unit crit}, we can tell which of these modules are unitary. Since $\La^{\dom}=\rho$ falls into case (2) of Corollary \ref{cor unit crit}, and the coordinates are integers, the unitary $\La$ are $\La^0=(-1,-2,\dots,-n)$ and 
\[
\La^q=(q,\dots,1,-q-1,\dots,-n),\qquad  q=1,2,\dots,n.
\]
Note that whenever $q\geq 1$, the inequality $\La^q_1\leq\frac{q+r}{2}$ is an equality, so $\La^q$ is the last point of unitarity on its line. On the other hand, $\La^0=(-1,\dots,-n)$ is in the continuous part of its line.

The Young diagram $Y^q$ corresponding to $\La^q$ has $n-q$ rows of length $n+1$ and then $q$ rows of length $n-q$. This includes the case $q=0$, with all rows of length $n+1$, and $q=n$ which corresponds to the empty diagram and the trivial module.
We can think of each $Y^q$ as a ``thick hook", or alternatively, as obtained by removing the $q\times (q+1)$ box from the lower right corner of the $n\times (n+1)$ box. 
In Example \ref{ex yd sp}, these diagrams are the full box, and the three diagrams on the diagonal consisting of lower left corners (i.e., the points such that the Hasse diagram contains no points to the left or below). We call the corresponding $\La^q$, except the one corresponding to the full box, the edge points of the Hasse diagram.
\medskip

We now consider general regular integral or half integral parameters, i.e., $\La^{\dom}$ that have integral or half integral coordinates with no repetitions and no zero.
By Corollary \ref{cor unit crit} (2), the unitary $\La$  belong to one of the cases (2a)-(2d). If $\La$ belongs to case (2b), then the negative coordinates must be $\leq -q-\half$; such $\La$ can be written as
\[
\La=(q-\half,\dots,\half,-q-\half-a_{n-q},-q-\frac{3}{2}-a_{n-q-1},\dots,-n+\half-a_1)
\]
for some integers $a_1\geq a_2\geq\dots\geq a_{n-q}\geq 0$.
These $\La$ belong to a reduced translation cone of dimension $n-q$ with vertex at 
\[
(q-\half,\dots,\half,-q-\half,-q-\frac{3}{2},\dots,-n+\half).
\]
If $\La$ belongs to case (2c), then the negative coordinates must be $\leq -q-1$; such $\La$ can be written as
\[
\La=(q,\dots,1,-q-1-a_{n-q},-q-2-a_{n-q-1},\dots,-n-a_1)
\]
for some integers $a_1\geq a_2\geq\dots\geq a_{n-q}\geq 0$.
These $\La$ belong to a reduced translation cone of dimension $n-q$ with vertex at the edge point
\[
(q,\dots,1,-q-1,-q-2,\dots,-n)
\]
of the Hasse diagram of $\rho$.

If $\La$ belongs to case (2d), then the negative coordinates after $-\half$ must be $\leq -q-\frac{3}{2}$; such $\La$ can be written as
\[
\La=(q+\half,\dots,\frac{3}{2},-\half,-q-\frac{3}{2}-a_{n-q-1},-q-\frac{5}{2}-a_{n-q-2},\dots,-n+\half-a_1)
\]
for some integers $a_1\geq a_2\geq\dots\geq a_{n-q-1}\geq 0$.
These $\La$ belong to a reduced translation cone of dimension $n-q-1$ with vertex at 
\[
(q+\half,\dots,\frac{3}{2},-\half,-q-\frac{3}{2},\dots,-n+\half).
\]
Finally, all $\La$ belonging to the case (2a) are of the form
\[
(-1-a_n,-2-a_{n-1},\dots,-n-a_1)\qquad\text{or}\qquad (-\half-a_n,-\frac{3}{2}-a_{n-1},\dots,-n+\half-a_1)
\]
for some integers $a_1\geq a_2\geq\dots\geq a_n\geq 0$. These $\La$ form two full translation cones, with vertices at
\[
(-1,-2,\dots,-n)\qquad\text{and}\qquad (-\half,-\frac{3}{2},\dots,-n+\half).
\]

\end{document}